\documentclass[11pt]{article}

\usepackage{amsmath} \allowdisplaybreaks
\usepackage{amssymb}
\usepackage{amsfonts}
\usepackage{dsfont}
\usepackage[table,xcdraw]{xcolor}
\usepackage{graphicx}
\usepackage{url}
\usepackage{booktabs}
\usepackage{multirow}
\usepackage{enumitem}
\usepackage{subcaption}
\usepackage{rotating}

\DeclareMathOperator*{\Minimize}{Minimize}
\DeclareMathOperator*{\Maximize}{Maximize}
\DeclareMathOperator*{\argmin}{argmin}

\begin{document}

\title{The Traveling Salesman Problem with non-convex polygonal neighborhoods: A hybrid Iterated Local Search–Block Coordinate Descent approach\thanks{This work has been partially supported by the Brazilian agencies FAPESP (grants 2022/05803-3 and 2023/08706-1), CNPq (grant 302073/2022-1) and CAPES (Finance Code 001).}}

\author{
    Ernesto G. Birgin\thanks{Department of Computer Science, Institute of Mathematics, Statistics and Computer Science, University of S\~ao Paulo, Rua do
    Mat\~ao, 1010, Cidade Universit\'aria, 05508-090, S\~ao Paulo, SP,
    Brazil. e-mail: egbirgin@ime.usp.br}
\and 
    Paula C. R. Ertel\thanks{Department of Applied Mathematics, Institute of Mathematics, Statistics and Computer Science, University of S\~ao Paulo, Rua do
    Mat\~ao, 1010, Cidade Universit\'aria, 05508-090, S\~ao Paulo, SP,
    Brazil. e-mail: paulaertel@ime.usp.br}
}

\date{November 24, 2026}

\maketitle

\begin{abstract}
The Traveling Salesman Problem with Neighborhoods (TSPN) generalizes the classical TSP by requiring a tour to visit a set of regions in the Euclidean plane rather than fixed point coordinates. It has applications in last-mile delivery, surveillance, and remote sensing in precision agriculture. In this paper, we address the Euclidean TSPN with arbitrary, non-convex polygonal neighborhoods---a highly complex setting due to the joint continuous and combinatorial optimization challenges. To solve the continuous subproblem of determining optimal visit points for a given visiting sequence, we propose a Block Coordinate Descent (BCD) method based on solving fast one-dimensional geometric subproblems. Computational experiments demonstrate that BCD computes optimal visit points for a fixed sequence in almost all instances in which optimality is certified, more than six orders of magnitude faster than a commercial solver solving the equivalent Mixed-Integer Quadratically Constrained Programming (MIQCP) formulation. This efficiency enables an effective hybrid Iterated Local Search metaheuristic (ILS-BCD) that thoroughly explores the combinatorial search space at low computational cost. Extensive computational tests on 235 benchmark instances derived from the Generalized TSP (GTSP), with up to 100 neighborhoods, show that ILS-BCD matched or improved upon the best-known GTSP tours in all cases, obtaining strictly shorter tours in 234 instances. On the 80 instances for which optimality is proven, ILS-BCD attains the optimal tour length in all of them, within a relative difference below $0.001\%$. On 60 larger instances, with as many as 6{,}325 neighborhoods, ILS-BCD remains competitive up to about 220.\\


\end{abstract}

\section{Introduction} \label{sec1}

In its classical form, the Traveling Salesman Problem (TSP) can be described as follows: given~\(m\) cities, with known travel distances between every pair of cities \(i\) and \(j\), the goal is to find a tour (i.e., an ordering of the cities) of minimum total length that visits each of the $m$ cities exactly once and returns to the starting city~\cite{Applegate2006}. When each city is represented by coordinates in the Euclidean plane and the distance between any two cities is defined as the Euclidean distance between their coordinates, the problem is known as the Euclidean TSP. Although simple to describe, the TSP belongs to the class of NP-hard problems~\cite{Garey1979} and is one of the most well-known problems in Combinatorial Optimization and Operations Research. Its prominence stems from its wide range of applications, the many variants it encompasses~\cite{GutinPunnen2007}, and the continued interest in the development of heuristic methods for solving it~\cite{Laporte2010}. Thus, the TSP has been, and continues to be, extensively studied. Nevertheless, some variants of the TSP remain relatively unexplored. One such variant is the Traveling Salesman Problem with Neighborhoods (TSPN).

In the Euclidean TSPN, each city is represented by a region in the Euclidean plane, referred to as a neighborhood or region. Finding a solution to this problem consists of two parts: determining (i) an ordering of the neighborhood and (ii) a visit point within each neighborhood. Now, these two parts in conjunction induce a tour that visits each city exactly once and returns to the starting city. The goal is to find (i) and (ii) so as to minimize the total length of the resulting tour. The TSPN also belongs to the class of NP-hard problems, since the TSP is the special case in which each region consists of a single point.

In the literature, the TSPN is sometimes referred to as the Close-Enough Traveling Salesman Problem (CETSP)~\cite{Gulczynski2006,Mennell2009}. However, the CETSP corresponds to the particular case in which the neighborhood of each city is a circle centered at its geographic coordinates. The problem can be interpreted as a Euclidean TSP in which the objective is not to visit each city exactly, but rather to travel within a predefined distance of each of them. In \cite{Behdani2014}, the authors formulate the CETSP as a mixed-integer nonlinear programming problem and propose methods to compute lower and upper bounds on the cost of an optimal tour. In \cite{Coutinho2016}, the authors propose a branch-and-bound algorithm that uses second-order cone programming to formulate the subproblem of determining the visit points within each neighborhood. Methods for the CETSP differ mainly in how they compute the visit points of a candidate sequence: by discretizing each disk into a finite set of candidates~\cite{Carrabs2017,Carrabs2017b}, by solving a second-order cone program for a fixed sequence, as in the Steiner-zone heuristic of~\cite{Mennell2011} and in the branch-and-bound of~\cite{Coutinho2016}, later improved in~\cite{Zhang2023}, or by exploiting the nonlinear structure of that subproblem to compute the points approximately~\cite{Wang2019}. Applications of the CETSP include radio-frequency identification technologies for automated meter reading~\cite{Dong2007,Shuttleworth2008}, last-mile delivery in which trucks need only come within the range of the drones they launch~\cite{Poikonen2017}, inspection of solar panels~\cite{DiPlacido2022} and data collection from sparsely deployed wireless sensors by a mobile robot~\cite{BoYuan2007}.

The work in \cite{Arkin1994} was the first to introduce the TSPN with arbitrary neighborhoods, under the name Geometric Covering Salesman Problem. In this study, the authors develop approximation algorithms for cases in which the neighborhoods have specific shapes, such as parallel line segments of equal length or neighborhoods defined as translations of a common convex polygon. These algorithms select a representative point for each neighborhood and then apply a classical TSP approximation algorithm. In \cite{Yang2018}, the authors proposed an algorithm for the TSPN that combines Particle Swarm Optimization to determine the visit points within each neighborhood with a genetic algorithm to generate the permutations. Although the method targets the TSPN with arbitrary neighborhoods, the authors show test results on instances that are all CETSP except for one. (This one instance contains neighborhoods of different shapes, including line segments, disks, rings, circular sectors, and some polygons.)

The papers \cite{Dror1999} and \cite{Hoeft1997} address the problem of optimizing the routing of industrial cutting machines for polygonal shapes on large plates, but from complementary perspectives. In both works, the $\ell_{\infty}$ norm is adopted for distance calculations, a metric that reflects the mechanics of real cutting machines, where two independent motion mechanisms operate simultaneously along the $x$ and $y$ axes. In this setting, the distance between two points is determined by the maximum displacement along a single coordinate. The paper \cite{Hoeft1997} introduces the Plate-Cutting Traveling Salesman Problem (P-TSP), an NP-hard problem focused on continuous cutting, in which the tool must enter and exit each polygon at exactly the same point. To tackle this problem, the authors propose heuristic methods based on Lagrangian decomposition, integrated with an iterative improvement algorithm, the variable r-opt heuristic, aimed at approximating high-quality routes and optimal cutting points. On the other hand, \cite{Dror1999} takes a complementary viewpoint: it focuses on the cutting subproblem in which the visitation order of the pieces is predetermined. This work advances previous research by proving that this subproblem (and its variants for convex polygons or entry and exit points restricted to vertices) can be solved optimally in polynomial time. To achieve this, the authors employ linear programming models and shortest-path algorithms on bipartite graphs, thereby eliminating the need for heuristic approaches at this stage. Extending this line of work, \cite{Dror2003} abstracts away specific industrial constraints (such as the use of the $\ell_{\infty}$ metric and the requirement that entry and exit points coincide) and generalizes the geometric problem by formulating the Touring Polygons Problem (TPP). The main objective of the TPP is to find the shortest Euclidean path that starts at a point $s$, visits a sequence of polygons in a predefined order, and ends at a destination $t$. The authors present the first exact algorithm capable of solving the problem in polynomial time for instances involving convex polygons. They further prove that when the polygons are non-convex, the problem becomes NP-hard. Since our problem involves non-convex geometries, the result of \cite{Dror2003} highlights the dual intractability of our setting: beyond the combinatorial explosion involved in determining the optimal visitation order, even computing the optimal entry and exit points for a fixed route is, in itself, NP-hard.

In this work, we propose an Iterated Local Search (ILS)~\cite{Loureco2010} algorithm to solve the Euclidean TSPN with polygonal neighborhoods. The algorithm combines eight classical TSP neighborhood searches, which generate permutations of a current solution, with a block coordinate descent (BCD) method~\cite{Amaral2022,Birgin2022,Wright2015} that determines the visit points of each permutation. The use of BCD for this purpose comes from~\cite{Birgin2022}, where it is illustrated by finding a short route that visits all neighborhoods of the city of São Paulo, Brazil, within an insertion-based constructive heuristic and a reinsertion local search. The most important difference lies in the subproblems of the BCD. The method of~\cite{Birgin2022} is a regularized BCD, whose subproblems have no known exact solution procedure; in the São Paulo application they are solved approximately, by discretization and brute-force search. Here, the subproblems are not regularized, and we describe how to compute their exact solution, analytically and through a numerically stable geometric formula. In addition, the BCD is embedded in a metaheuristic with a much richer set of neighborhoods, and the resulting method is evaluated on a large set of benchmark instances. We chose ILS motivated by the robust results reported in \cite{Schmidt2022} for the Generalized Traveling Salesman Problem (GTSP)~\cite{Fischetti1997}. In a companion paper~\cite{Birgin2026cetsp}, we consider this same combination for the CETSP, the case in which every neighborhood is a disk. There, the BCD subproblems are solved by a different method, specialized for disks, and the BCD is also combined with the population-based Lin--Kernighan--Helsgaun heuristic.

The remainder of this work is organized as follows. In Section~\ref{sec2}, we formally define the problem by presenting a mixed-integer quadratically constrained programming formulation. Section~\ref{sec3} explains how to use the BCD method to compute the visit points for a given permutation of the cities. Section~\ref{sec4} presents the ILS method and the eight TSP neighborhoods used to generate neighboring permutations of a current solution. In the experiments section, we report numerical results on a large set of adapted GTSP instances from the literature. We first describe how we adapt the GTSP instances for the TSPN. Next, we demonstrate the efficiency of the BCD method for computing visit points of a given permutation, comparing it with CPLEX on the same task. Beyond efficiency, we analyze the BCD method’s ability to compute optimal visit points and highlight cases where it may find only local solutions. Finally, we present extensive experiments solving TSPN instances with non-convex polygonal neighborhoods and compare the obtained solutions with known solutions of the GTSP instances associated with each TSPN instance. On the subset of instances for which the MIQCP formulation, solved with CPLEX, proves optimality, we also compare the tours returned by the metaheuristic with the optimal ones. We close the section by examining how the method behaves on substantially larger instances and where it ceases to be competitive. Section~\ref{sec6} summarizes the contributions and indicates directions for future research.\\

\noindent
\textbf{Notation.} We denote by $\mathbb{R}^n$ the $n$-dimensional Euclidean space, and by $\|\cdot\|$ the Euclidean norm, with $\langle \cdot, \cdot \rangle$ representing the Euclidean inner product. For a vector $x$, $x^T$ is its transpose. For a set $S$, $\mathrm{int}(S)$ and $\partial S$ denote its interior and boundary, respectively. The line segment connecting points $p$ and $q$ is denoted by $\overline{pq}$. Set operations are represented by $\cap$ (intersection) and $\cup$ (union). The empty set is represented by $\emptyset$.

\section{Problem description and model} \label{sec2}

\begin{figure}[ht!]
\begin{center}
\begin{tabular}{ccc}
\includegraphics{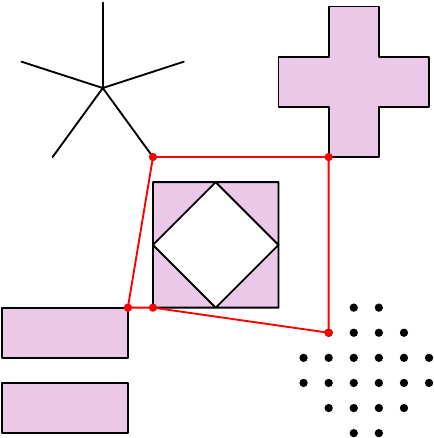} & &
\includegraphics{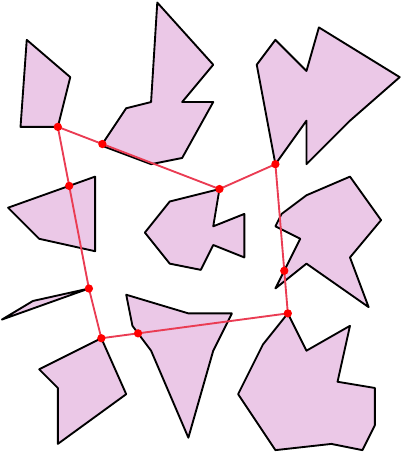} \\ [2mm]
(a) Five-city instance & & (b) Ten-city instance
\end{tabular}
\end{center}
\caption{In this study, a neighborhood is defined as the union of a finite number of convex polygons. Figures~(a) and~(b) show two groups of such neighborhoods. Figure (a) illustrates the level of generality that the neighborhoods' shape can take. It is a five-city instance in which a neighborhood may have no interior, be the union of a finite number of points, be disconnected, or be a non-convex or non-simple polygon. Figure (b) illustrates a ten-city instance in which, for each neighborhood, the union of the convex polygons constitutes a non-convex simple polygon. In both cases, the optimal tour is displayed.} 
\label{fig1}
\end{figure}

In this section we present a mixed-integer quadratically constrained programming (MIQCP) formulation for the TSPN with polygonal regions and Euclidean distances. An instance of the TSPN is defined by a set of input data as follows. The number of cities is given by $m$. Because each---generally non-convex---polygonal neighborhood $\mathcal{N}_i$ of city~$i$ can be expressed as the union of a set of convex polygons, we have that $\mathcal{N}_i = \cup_{\ell=1}^{m_i} Q_{i\ell}$, where each $Q_{i\ell}$ is a convex polygon with $q_{i\ell}$ vertices $v_{i\ell}^1,\dots,v_{i\ell}^{q_{i\ell}}$. Thus each $Q_{i\ell}$ with $q_{i\ell} \geq 3$ can be easily described as $Q_{i\ell} = \{ (x,y) \in \mathbb{R}^2 \;|\; a_{i\ell k} \, x + b_{i\ell k} \, y + c_{i\ell k} \leq 0, k = 1,\dots, q_{i\ell} \}$, where the coefficients of the linear inequalities depend on its vertices. If $Q_{i\ell}$ is a single point, then $Q_{i\ell} = \{ v_{i\ell}^1 \}$; and if $Q_{i\ell}$ is a segment, then $Q_{i\ell} = \{ x\in \mathbb{R}^2 \;|\; x =  v_{i\ell}^1 + \lambda_{i\ell} (v_{i\ell}^2 - v_{i\ell}^1), \lambda_{i\ell} \in [0,1] \}$. As defined, each neighborhood is composed by the union of a finite number of convex polygons; see Figure~\ref{fig1}. Considering this definition of neighborhood, the TSPN can be viewed as a generalization of the GTSP~\cite{Fischetti1997}, in which each city is represented by a finite set of points. In turn, the GTSP generalizes the classical TSP, in which each city is represented by a single point.

The formulation involves several types of decision variables. Each neighborhood \(\mathcal{N}_i\) has an associated visit point \(p_i = (p_i^x, p_i^y)^T \in \mathbb{R}^2\) whose coordinates are decision variables representing the point to be visited within that neighborhood. The binary variables \(w_{i\ell}\) are used to select exactly one polygon \(Q_{i\ell}\) within which the visit point \(p_i\) must lie. When the polygon \(Q_{i\ell}\) is defined as a line segment, the continuous variable \(\lambda_{i\ell} \in [0,1]\) determines the position of \(p_i\) along the segment. The binary variables \(x_{ij}\) indicate the routing decisions, where \(x_{ij} = 1\) if the tour travels directly from \(p_i\) to \(p_j\), and \(x_{ij} = 0\) otherwise. The continuous variables \(u_i\) are sequencing variables that help define the visiting order of the points and ensure the construction of a single tour. The variables \(\xi_{ij}\) represent the Euclidean distances between the visit points \(p_i\) and \(p_j\), while \(\gamma_{ij}\) denote the corresponding travel costs associated with the selected arcs of the tour.  Together, these variables determine both the spatial configuration of the selected visit points and the optimal route that connects them. The model follows below.

{\footnotesize
\begin{align}
\Minimize & 
\sum_{i=1}^m \sum_{\substack{j=1\\j\neq i}}^m \gamma_{ij} \label{g1}\\
\text{s.t. } 
& \sum_{\substack{j=1\\j\neq i}}^{m} x_{ij} =1 & i = 1,\dots,m \label{g2}\\
&\sum_{\substack{i=1\\i\neq j}}^{m} x_{ij} =1 & j = 1,\dots,m \label{g3}\\
& u_i - u_j + (m-1) x_{ij} \leq m-2 & i,j=2,\dots,m, i\neq j \label{g4}\\
& u_i \geq 1 & i=2,\dots,m \label{g5}\\
& u_i \leq m-1 & i=2,\dots,m \label{g6}\\
& \xi_{ij}^2 \geq (p_i^x - p_j^x)^2 + ( p_i^y - p_j^y)^2 & i,j=1,\dots,m, i\neq j \label{g7}\\
& \xi_{ij} \geq 0 & i,j=1,\dots,m, i\neq j \label{g8}\\
& \gamma_{ij} \leq M x_{ij} & i,j=1,\dots,m, i\neq j \label{g9}\\
& \gamma_{ij} \geq \xi_{ij} - M (1-x_{ij}) & i,j=1,\dots,m, i\neq j \label{g10}\\
& \gamma_{ij} \geq 0 & i,j=1,\dots,m, i\neq j \label{g11}\\
& p_i \geq v_{i\ell}^{1} - M(1-w_{i\ell}) \, e & i = 1,\dots,m, \ell = 1,\dots,m_i \;|\; q_{i\ell}=1 \label{g12}\\
& p_i \leq v_{i\ell}^{1} + M(1-w_{i\ell}) \, e & i = 1,\dots,m, \ell = 1,\dots,m_i \;|\; q_{i\ell}=1 \label{g13}\\
& p_i \geq v_{i\ell}^{1} + \lambda_{i\ell} (v_{i\ell}^{2} - v_{i\ell}^{1}) - M(1-w_{i\ell}) \, e & i = 1,\dots,m, \ell = 1,\dots,m_i \;|\; q_{i\ell}=2 \label{g14}\\
& p_i \leq v_{i\ell}^{1} + \lambda_{i\ell} (v_{i\ell}^{2} - v_{i\ell}^{1}) + M(1-w_{i\ell}) \, e & i = 1,\dots,m, \ell = 1,\dots,m_i \;|\; q_{i\ell}=2 \label{g15}\\
& a_{i\ell k} p_i^x + b_{i\ell k} p_i^y + c_{i\ell k} \leq M(1-w_{i\ell}) & i=1,\dots,m, \ell = 1,\dots,m_i, k = 1,\dots, q_{i\ell} \;|\; q_{i\ell} \geq 3 \label{g16}\\
& \sum_{\ell=1}^{m_i} w_{i \ell} = 1 & i = 1,\dots,m \label{g17}\\
& x_{ij} \in \{0,1\} & i,j=1,\dots,m, i\neq j  \label{g18}\\
& u_i \in \mathds{R} & i=2,\dots,m, \label{g19}\\
& \xi_{ij} \in \mathds{R} & i,j=1,\dots,m, i\neq j \label{g20}\\
& \gamma_{ij} \in \mathds{R} & i,j=1,\dots,m, i\neq j \label{g21}\\
& \lambda_{i\ell} \in [0,1] & i = 1,\dots,m, \ell = 1,\dots,m_i \;|\; q_{i \ell} = 2 \label{g22}\\
& w_{i\ell} \in \{0,1\} & i = 1,\dots,m, \ell = 1,\dots,m_i \label{g23}\\
& p_i \equiv (p_i^x, p_i^y)^T \in \mathbb{R}^2 & i = 1,\dots,m \label{g24}
\end{align}}

The constraints in the model can be grouped according to their roles. Constraints~\eqref{g2} and~\eqref{g3} ensure that the tour visits each neighborhood exactly once by requiring exactly one outgoing and one incoming arc per visit point. The subtour elimination is enforced by constraints~\eqref{g4}, together with the bounds on the sequencing variables \(u_i\) in~\eqref{g5} and~\eqref{g6}; these constraints implement a Miller--Tucker--Zemlin~\cite{MTZ1960} formulation to guarantee a valid visiting order and prevent subtours. Constraints~\eqref{g7} and~\eqref{g8} define the Euclidean distance \(\xi_{ij}\) between the visit points \(p_i\) and \(p_j\), while constraints~\eqref{g9}--\eqref{g11} link these distances to the selected arcs through the variables~\(\gamma_{ij}\). The placement of visit points within their neighborhoods is ensured by constraints~\eqref{g12}--\eqref{g17}, handling cases of fixed points, line segments, polygons, and multipolygonal regions. Finally, constraints~\eqref{g18}--\eqref{g24} define the domains of all decision variables, specifying which variables are binary and which are continuous, as well as the bounds imposed on the latter. The objective function of the model, given by~\eqref{g1}, aims to minimize the total travel distance of the tour by summing the effective distances \(\gamma_{ij}\) over all selected arcs. Because \(\gamma_{ij}\) is constrained to be positive only when \(x_{ij} = 1\) through constraints~\eqref{g9}--\eqref{g11}, the objective effectively measures the total length of the tour that visits each neighborhood exactly once. In summary, the model determines (i) the optimal order in which to visit a set of neighborhoods, (ii) the precise visit point within each neighborhood, and (iii) the total minimum travel distance for the resulting tour. It represents therefore a traveling salesman problem with (non-convex) polygonal neighborhoods, integrating both combinatorial routing and continuous spatial optimization in a single framework.

In constraints~(\ref{g9},\ref{g10}) and~\eqref{g12}--\eqref{g16}, $M>0$ is a constant that must assume a sufficiently large value. All neighborhoods $\mathcal{N}_i$ are contained within a rectangle, say, $[\underline{p}^x, \overline{p}^x] \times [\underline{p}^y, \overline{p}^y]$. In~(\ref{g9},\ref{g10}), we can define $M$ as $\Vert (\overline{p}^x,\overline{p}^y) - (\underline{p}^x,\underline{p}^y) \Vert$, which is an upper bound on the distance between any two visit points. In \eqref{g16}, it would suffice to set $M$ as the optimal value of
\begin{equation*}
    \Maximize_{i=1,\dots,m, \ell = 1,\dots,m_i, k = 1,\dots, q_{i \ell}} \left\{ \Maximize_{\{ (p^x,p^y)^T \in [\underline{p}^x, \overline{p}^x] \times [\underline{p}^y, \overline{p}^y] \}} \{ a_{i\ell k}  p_i^x + b_{i\ell k}  p_i^y + c_{i\ell k} \} \right\}.
\end{equation*}

Based on the well-known norm inequality $\Vert x \Vert_\infty \leq \Vert x \Vert_2 \leq \Vert x \Vert_1$ for $x \in \mathds{R}^n$, bounds on the distances between pairs of visit points can be incorporated into model~(\ref{g1})--(\ref{g24}). Through the introduction of additional continuous variables $\xi^{x}_{ij}$ and $ \xi^{y}_{ij}$ for $i,j=1,\dots,m, i\neq j$, the following inequalities can be imposed:
\begin{align}
     - \xi^{x}_{ij} \leq p^{x}_{i} - p^{x}_{j} & \leq \xi^{x}_{ij} & i,j=1,\dots,m, i\neq j \label{gac1}\\
     - \xi^{y}_{ij} \leq p^{y}_{i} - p^{y}_{j} & \leq \xi^{y}_{ij} & i,j=1,\dots,m, i\neq j \label{gac2}\\
     \xi^{x}_{ij}, \xi^{y}_{ij} & \geq 0 & i,j=1,\dots,m, i\neq j \label{gac3}\\
     \xi_{ij} & \geq \xi^{x}_{ij} & i,j=1,\dots,m, i\neq j \label{gac4}\\
     \xi_{ij} & \geq \xi^{y}_{ij} & i,j=1,\dots,m, i\neq j \label{gac5}\\
     \xi_{ij} \leq \xi^{x}_{ij} & + \xi^{y}_{ij} & i,j=1,\dots,m, i\neq j. \label{gac6}
\end{align}
Constraints~(\ref{gac4},\ref{gac5}) correspond to \( \| p_i - p_j \|_2 \geq \| p_i - p_j \|_\infty \) while constraints~\eqref{gac6} correspond to \( \| p_i - p_j \|_2 \leq \| p_i - p_j \|_1 \). Including these inequalities adds a total of $2m(m-1)$ nonnegative continuous variables and $7m(m-1)$ linear constraints to model~\eqref{g1}--\eqref{g24}.

\section{Block coordinate descent} \label{sec3}

In this section, we consider the problem of, given a permutation $i_1,\dots,i_m$ of the cities, determining the visit points $p_{i_j} \in \mathcal{N}_{i_j}$ for $j=1,\dots,m$ that minimize the length of the tour starting at $p_{i_1}$, visiting $p_{i_2},\dots,p_{i_m}$ in this order, and returning to $p_{i_1}$. Specifically, we describe how to solve this problem using the block coordinate descent (BCD) method. For the remainder of this section, we simplify the notation by defining $i_0 = i_m$ and $i_{m+1}=i_1$. To simplify the exposition, we assume that each neighborhood is a single point, a segment, or a (possibly non-convex) simple polygon, as depicted in Figure~\ref{fig1}(b). The strategy we describe here, however, can be trivially adapted to all the neighborhood types introduced in Section~\ref{sec2} and illustrated in Figures~\ref{fig1}(a) and~\ref{fig1}(b).

The Block Coordinate Descent (BCD) method is a generalization of the standard coordinate descent method~\cite{Wright2015}. In BCD, instead of updating a single variable at a time, blocks of variables are updated simultaneously. At each iteration, a block of variables, either predefined or chosen adaptively, is selected, and a function (which does not necessarily coincide with the original objective function) is minimized with respect to that block while keeping the remaining variables fixed. This process is repeated for successive blocks until convergence. The method is particularly useful when each block’s subproblem is easy to solve, enabling large or complex optimization problems to be tackled efficiently by breaking them into simpler, more manageable pieces. In practice, the algorithm starts with an initial guess, iteratively selects a block of variables, optimizes over this block of variables, and continues cycling through or adaptively choosing blocks until the solution stabilizes or the change in the objective function becomes sufficiently small. The works~\cite{Amaral2022} and~\cite{Birgin2022} introduce BCD methods for minimizing smooth nonconvex functions subject to bound constraints and general constraints, respectively. The main feature of these methods is that, at each iteration, they minimize a regularized function. Both works present complexity results in terms of the number of iterations and function evaluations required to find approximate stationary points.

The visit points that minimize the tour length for the permutation $i_1,\dots,i_m$ are a solution to the problem
\begin{equation} \label{bcdprob}
\Minimize_{p_{i_j} \in \mathcal{N}_{i_j},\, j=1,\dots,m} \Psi(p_{i_1},\dots,p_{i_m}) := \sum_{j=1}^m \| p_{i_j} - p_{i_{j-1}} \|.
\end{equation}
Problem~\eqref{bcdprob} is particularly well-suited for a BCD method. By associating each visit point with a block, the optimization can be performed block by block, updating the visit point of one city at a time while keeping the others fixed. In this way, each subproblem reduces to updating a single point $p_{i_{\ell}} \in \mathbb{R}^2$ as
\begin{equation} \label{bcdsubprob}
p_{i_{\ell}} = \argmin_{\pi \in \mathcal{N}_{i_{\ell}}} \left\{ \phi(\pi) := \| p_{i_{\ell+1}} - \pi \| + \| \pi - p_{i_{\ell-1}} \| \right\}.
\end{equation}
Since all the other visit points are fixed, minimizing the objective function in~\eqref{bcdsubprob} is equivalent to minimizing the objective function in~\eqref{bcdprob}. The method of~\cite{Birgin2022} accepts the update of a block only if it satisfies a sufficient descent condition, and otherwise solves a regularized subproblem; this safeguard, together with the regularization also used in~\cite{Amaral2022}, is what yields their guarantees of convergence to stationary points. The method defined here computes the global minimizer of~\eqref{bcdsubprob} without checking for sufficient descent or adding a regularization term, so the convergence analysis of~\cite{Birgin2022} does not apply to it. In return, as we show next, problem~\eqref{bcdsubprob} can be solved exactly, in closed form on each edge of the neighborhood, whereas the subproblems of~\cite{Birgin2022} have no known exact solution procedure and, in its São Paulo application, are solved by discretization and brute-force search.

Two hypotheses on the polygonal neighborhoods are enough to turn~\eqref{bcdprob} into a problem on which the BCD method just described is guaranteed to find a global minimizer. Suppose first that every $\mathcal{N}_{i_j}$ is convex. The feasible set $\mathcal{N} = \mathcal{N}_{i_1} \times \dots \times \mathcal{N}_{i_m}$ is then a convex and compact subset of $\mathbb{R}^{2m}$, and $\Psi$ in~\eqref{bcdprob} is convex and continuous on it. Suppose in addition that consecutive polygons in the permutation do not meet. Thus, no two consecutive visit points can coincide and $\Psi$ is continuously differentiable on $\mathcal{N}$. A convex function that is differentiable is also pseudoconvex, which is what the argument below needs. Under these two hypotheses, solving each subproblem~\eqref{bcdsubprob} to global optimality makes every update well defined and never increases~$\Psi$, so the iterates form an infinite sequence with non-increasing objective values, confined to the compact set $\mathcal{N}$. In~\cite{Grippo2000convergence}, the authors proved that, for a pseudoconvex objective on a Cartesian product of convex sets, the limit points of such a sequence are global minimizers, and the visit points delivered by the method would therefore be optimal for the given permutation. 

We now turn our attention to solving problem~\eqref{bcdsubprob}. To simplify notation, let \(p \equiv p_{i_{\ell+1}}\), \(q \equiv p_{i_{\ell-1}}\), and \(P \equiv \mathcal{N}_{i_{\ell}}\). Then, the optimization problem in~\eqref{bcdsubprob} can be rewritten as
\begin{equation} \label{bcdsubprob3}
\Minimize_{\pi \in P} \phi(\pi) := \| p - \pi \| + \| \pi - q \|.
\end{equation}
If \(\overline{pq} \cap P \neq \emptyset\), then any \(\bar \pi \in \overline{pq} \cap P\) is optimal, since the triangle inequality implies \(\phi(\pi) \ge \|p-q\| = \phi(\bar \pi)\) for all \(\pi \in \mathbb{R}^2\). Hence, we only need to consider the case \(\overline{pq} \cap P = \emptyset\). In this case, no \(\pi \in \mathrm{int}(P)\) can be optimal: \(\phi\) is differentiable at every \(\pi \notin \{p,q\}\), and its gradient \((\pi-p)/\|\pi-p\| + (\pi-q)/\|\pi-q\|\) vanishes only on \(\overline{pq}\), which does not meet \(P\). The search can therefore be restricted to \(\pi \in \partial P\). Let \(h\) denote the number of vertices of \(P\), labeled \(v_1, \dots, v_h\), and define \(v_0 = v_h\) and \(v_{h+1} = v_1\) for convenience. A possible parameterization of the boundary is then given by \(\pi(\lambda) = v_1 + (\lambda-1)(v_2 - v_1)\) with \(\lambda \in [1,h]\) if \(h \le 2\), and by \(\pi(\lambda) = v_i + (\lambda-i)(v_{i+1} - v_i)\) with \(\lambda \in [i,i+1]\) and $i=1,\dots,h$ if \(h \ge 3\). With these parameterizations, problem~\eqref{bcdsubprob3} can be rewritten as the univariate problem
\begin{equation} \label{bcdsubprob4}
\Minimize_{\lambda \in \mathbb{R}} \psi(\lambda) := \| p - \pi(\lambda) \| + \| \pi(\lambda) - q \| \mbox{ subject to } 1 \leq \lambda \leq \bar h,
\end{equation}
where $\bar h = h$ if $h \leq 2$ and $\bar h = h+1$ if $h \geq 3$.

The feasible region \(\{\lambda \in \mathbb{R} \mid 1 \le \lambda \le \bar h\}\) of~\eqref{bcdsubprob4} can be decomposed into the points \(\lambda = i\) for \(i = 1, \dots, \bar h\) and the intervals \(\lambda \in (i,i+1)\) for \(i = 1, \dots, \bar h-1\). The former correspond to the case where \(\pi(\lambda)\) coincides with the vertex \(v_i\), while the latter correspond to the case where \(\pi(\lambda)\) lies in the interior of the edge \((v_i, v_{i+1})\). The solution method consists of, for each edge, checking whether there exists a point \(\pi(\lambda)\) in its interior such that \(\psi'(\lambda) = 0\). The overall solution is then given by the point that attains the smallest value of \(\psi\) among the vertices and the interior points of the edges where \(\psi'\) vanishes. It remains only to show how to determine whether an edge contains a point $\pi(\lambda)$ where $\psi'(\lambda) = 0$.

For \(i=1,\dots,\bar h - 1\), if \(\lambda \in (i,i+1)\), the derivative of \(\psi\) is
\[
\psi'(\lambda) = \Big\langle - \frac{p-\pi(\lambda)}{\|p-\pi(\lambda)\|} + \frac{\pi(\lambda)-q}{\|\pi(\lambda)-q\|},\, v_{i+1}-v_i \Big\rangle,
\]
and therefore \(\psi'(\lambda)=0\) if and only if
\begin{equation} \label{nullpsider2}
\frac{\langle \pi(\lambda)-q, v_{i+1}-v_i \rangle}{\|\pi(\lambda)-q\|} = 
\frac{\langle p-\pi(\lambda), v_{i+1}-v_i \rangle}{\|p-\pi(\lambda)\|}.
\end{equation}
By squaring both sides (to eliminate the square roots) and performing the cross product, one obtains a quartic equation in~\(\lambda\). In fact, the quartic has its cubic and quartic coefficients equal to zero, reducing it to a quadratic equation whose solutions can be easily computed. However, the computation is not entirely straightforward, as several cases must be carefully considered to avoid numerically unstable calculations. We now show how using a geometrical perspective simplifies the solution of~\eqref{nullpsider2}. We start by noting that~\eqref{nullpsider2} is equivalent to
\(\cos(\theta_1) = \cos(\theta_2)\), where \(\theta_1\) is the angle between \(p - \pi(\lambda)\) and \(d = v_{i+1} - v_i\), and \(\theta_2\) is the angle between \(\pi(\lambda) - q\) and \(d\). Define the line \(L = \{ v_i + \tilde\lambda\, d \mid \tilde\lambda \in \mathbb{R} \}\), and the scalar projection of a point \(x\) along \(L\) as \(\alpha(x) = \langle x - v_i, d \rangle / \|d\|^2\). Then the perpendicular distances from \(p\) and \(q\) to \(L\) are \(\|p_\perp\| = \| p - v_i - \alpha(p) \, d \|\) and \(\|q_\perp\| = \| q - v_i - \alpha(q) \, d \|\), respectively. With these definitions, equation~\eqref{nullpsider2} reduces to the simple one-dimensional relation
\[
\frac{\tilde\lambda - \alpha(q)}{\alpha(p) - \tilde\lambda} = \frac{\|q_\perp\|}{\|p_\perp\|},
\]

whose solution is
\[
\tilde\lambda^* = \frac{\alpha(q) \|p_\perp\| + \alpha(p) \|q_\perp\|}{\|p_\perp\| + \|q_\perp\|}.
\]

Thus, the search for \(\lambda\) such that \(\psi'(\lambda)=0\), which originally appeared to require solving a quartic/quadratic equation, reduces to a one-dimensional geometric problem whose unique solution corresponds to a weighted interpolation of the projected points, with weights given by the perpendicular distances. If \(\tilde\lambda^* \in (0,1)\), the point \(v_i + \tilde\lambda^* d = \pi(i+\tilde\lambda^*)\) lies in the interior of the edge and is a candidate solution to~\eqref{bcdsubprob4}; otherwise \(\psi'\) does not vanish on the edge, and only its endpoints remain as candidates. The formula holds whether \(p\) and \(q\) lie on the same side of \(L\), in which case \(v_i + \tilde\lambda^* d\) is the point where \(L\) meets the segment joining \(p\) to the reflection of \(q\) across \(L\), or on opposite sides, in which case it is the point where \(\overline{pq}\) crosses \(L\); since \(\overline{pq} \cap P = \emptyset\), this crossing never lies on the edge itself. It requires \(\|p_\perp\| + \|q_\perp\| > 0\): if \(p\) and \(q\) both lie on \(L\), they lie on the same side of the edge, again because \(\overline{pq} \cap P = \emptyset\), so \(\psi\) is monotone on the edge and only its endpoints are candidates. 

\subsection{Discussion on block definition and convergence to a global minimizer}

In this section, we discuss several aspects that affect the performance of the BCD method when applied to problem~\eqref{bcdprob}. These aspects include the block definitions, the non-convexity of the cities' neighborhoods, and the non-differentiability of the objective function of~\eqref{bcdprob} when two consecutive visit points coincide.

\begin{figure}[ht!]
\begin{center}
\begin{tabular}{ccc}
\includegraphics{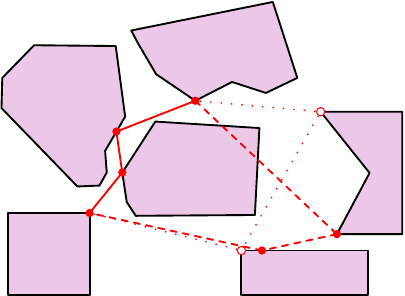} & &
\includegraphics{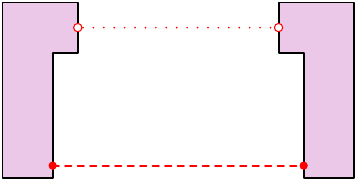} \\ [2mm]
\end{tabular}
\end{center}
\caption{This figure illustrates tours that cannot be improved when each block is associated with a single visit point, but can be improved when each block is instead associated with two consecutive visit points. The tours that cannot be improved under the single-visit-point assignment correspond to the arcs shown with solid and dashed lines. The improved tours obtained by associating each block with two consecutive visit points correspond to the arcs shown with solid and dotted lines. Notably, the tours differ only at two consecutive visit points. The repositioned visit points are shown using hollow markers.}
\label{fig2}
\end{figure}

Assigning each block to a single visit point simplifies the subproblem addressed by the BCD method at each iteration; however, this choice is not without drawbacks. Figure~\ref{fig2} illustrates two tours that cannot be improved under this block assignment scheme. However, these tours can be improved if we instead associate each block with two consecutive visit points. In this case, the subproblem to be solved at each iteration of the BCD method is given by
\begin{equation} \label{bcdsubprob2}
(p_{i_{\ell}}, p_{i_{\ell+1}}) = \argmin_{\hat \pi \in \mathcal{N}_{i_{\ell}}, \, \tilde \pi \in \mathcal{N}_{i_{\ell+1}}} \left\{ \| p_{i_{\ell+2}} - \tilde \pi \| + \| \tilde \pi - \hat \pi  \| + \| \hat \pi - p_{i_{\ell-1}} \| \right\}.
\end{equation}
This problem involves four variables instead of just two, is geometrically more challenging, and we are not aware of any method to solve it that is as efficient as the one we described for solving~\eqref{bcdsubprob}. To assess whether developing a dedicated method to solve~\eqref{bcdsubprob2} would be worthwhile, we conducted preliminary experiments. These tests involved discretizing the boundaries of both polygons by placing points along them and evaluating all point combinations between the two polygons. Although this brute-force strategy is computationally expensive, it revealed that, within the context of solving the TSPN, the gains in solution quality obtained by considering subproblems of the form~\eqref{bcdsubprob2} are only marginal. Furthermore, this brute-force strategy relies on discretization, an approach we seek to avoid. Developing an efficient method for the two-visit-point subproblem~\eqref{bcdsubprob2} that exploits the continuity of the feasible space remains a topic for future research.

The convergence of the BCD method to local solutions, such as those illustrated in Figure~\ref{fig2}, strongly depends on the initial guess of the iterative optimization process. A classical technique to increase the likelihood of converging to a global minimizer is the use of a multi-start strategy—that is, running the algorithm from multiple initial guesses. The application of this technique will be empirically analyzed in Section~\ref{sec532}. Still regarding local versus global minimizers, another alternative is applying BCD to subproblems with convex polygons. Up to this point, BCD has been presented as computing visit points given a fixed visiting sequence. In Section~\ref{sec2}, when introducing the MIQCP model for the TSPN, we assumed each non-convex polygonal neighborhood is the union of a finite set of convex subpolygons. However, this partition was not incorporated into our solution algorithm. Had we included it, a discrete optimization module could provide BCD with both a visiting sequence and the designated convex subpolygon for each city. BCD would then only determine the visit point inside that predefined subpolygon. Although this would prevent BCD from converging to local minimizers, it would expand the discrete search space for the combinatorial solver. We leave this strategy for future work on heuristic and exact methods for the TSPN.

The final issue to be addressed here is the non-differentiability of the objective function in~\eqref{bcdprob}, which occurs whenever two consecutive visit points in the sequence coincide. If the cities' polygonal neighborhoods are disjoint, the objective function remains differentiable at all BCD iterates because the method always generates feasible points. Similarly, even if neighborhoods overlap, the objective function remains differentiable across all iterates as long as overlapping polygons never appear consecutively in the visiting sequence. Finally, when two consecutively visited neighborhoods have a non-empty intersection, their visit points may or may not coincide: if they do not, differentiability is preserved; if they do coincide (or become numerically extremely close), the objective function becomes non-differentiable. We empirically analyze what occurs in these cases in Section~\ref{sec531}.

\section{Metaheuristic, local search and neighborhood structures} \label{sec4}

To solve the TSPN problem, we implemented an Iterated Local Search (ILS) metaheuristic~\cite{Loureco2010} in which the visit points of the candidate solutions are computed with the BCD method of Section~\ref{sec3}. We refer to the resulting hybrid metaheuristic as ILS-BCD. A solution is represented as a permutation of the cities $s = (s_1, \dots, s_m)$ together with the corresponding visit points $w = (w_{s_1}, \dots, w_{s_m})$ within their neighborhoods. The ILS begins with a solution obtained from a constructive heuristic and then performs a local search from this solution. The resulting solution $(s, w)$ then becomes the current solution. Each ILS iteration comprises the following sequential steps: perturbation of the current solution, a local search from the perturbed solution, and acceptance of the local search output $(\bar{s}, \bar{w})$ as the new current solution if a criterion is satisfied. A solution with a lower cost than the current solution $(s, w)$ is always accepted. Solutions with higher cost may also be accepted if
\[
f(\bar{s}, \bar{w}) \le f(s^*, w^*) + \eta |f(s^*, w^*)|,
\]
where $(s^*, w^*)$ is the best known solution and $f(s,w)$ represents the cost of a solution $(s,w)$. The threshold $\eta$ starts at $\eta_0$ and is multiplied by $h \in (0,1]$ after every $n_{\mathrm{stag}}$ consecutive iterations that fail to improve the best solution found. Once $\eta$ falls below $\eta_{\min} = 10^{-4}$ the acceptance test admits essentially nothing but improvements, and the search is reheated: $\eta$ returns to $\eta_0$ and the current solution is reset to the best one found so far, a restart acceptance criterion in the sense of~\cite{Loureco2010}. The parameters $\eta_0 > 0$, $h \in (0,1]$ and $n_{\mathrm{stag}} \geq 1$ require calibration, whereas $\eta_{\min}$ is fixed. With $h = 1$ the threshold never decreases and the reheating is never activated.

The initial solution is generated using the Nearest Neighbor (NN) heuristic~\cite{Rosenkrantz1977}, applied to a TSP in which each city is represented by the centroid of its polygonal neighborhood---recall that in the TSPN, cities are represented by neighborhoods rather than by their geographical coordinates. After determining the permutation using the NN heuristic, the visit points are recomputed using the BCD method to obtain a tour. By initializing the NN algorithm from each city, and applying the BCD to each of the $m$ tours built this way, the best resulting tour and its corresponding visit points are selected as the starting solution.

For the local search step, we use a Variable Neighborhood Descent (VND) metaheuristic~\cite{Hansen2010}, which explores the solution space using multiple neighborhood structures. Given a set of $\kappa$ neighborhoods $\mathcal{N}_1, \dots, \mathcal{N}_{\kappa}$ and a current solution $(s, w)$, a VND iteration begins with $\mathcal{N}_1$. If a neighbor $(\bar{s}, \bar{w})$ that improves upon the current solution is found, it becomes the new current solution and the iteration ends. Otherwise, the next neighborhood, $\mathcal{N}_2$, is explored, continuing sequentially up to $\mathcal{N}_{\kappa}$. Whenever the current solution changes, the search restarts with $\mathcal{N}_1$. In our implementation we adopt the $\kappa = 8$ neighborhoods described below, and their exploration order in the VND is subject to calibration: fixed or drawn at random at each call. Moreover, within a neighborhood, whether the search returns the first improving neighbor or the best one is also calibrated for the 2-Opt, 3-Opt, swap, insertion and string relocation neighborhoods. This choice does not apply to the double-bridge, Balas--Simonetti and Gutin neighborhoods, which are not explored exhaustively: for each of them a dedicated procedure finds the best neighbor for the current visit points, which are kept fixed during the search, and only that neighbor is then evaluated with the BCD.

The solution neighborhoods listed below were adapted for the TSPN and implemented in the VND. Each description corresponds to a type of move. By applying a single move to the current solution $(s, w)$, we obtain a permutation that we call a neighbor of $(s, w)$. The set of all possible neighbors---each obtained by applying a single move of the same type to $(s, w)$---form a neighborhood of the current solution. As we have just described, a move generates a new permutation by modifying the current one, but this constitutes only the first stage of neighbor construction. In the second stage, all visit points—not only those associated with the cities that have been moved—must be recalculated. The neighbor’s tour length can be computed only after both stages are complete. Evaluating a neighbor (i.e., computing its tour length) requires recalculating all visit points $p$, not only those corresponding to the cities that were moved, because even a single change in the permutation can affect the optimal location of every visit point. In the 2-Opt, 3-Opt, swap, insertion and string relocation neighborhoods, every neighbor examined is evaluated in this way, with all its visit points recalculated by the BCD. In the double-bridge, Balas--Simonetti and Gutin neighborhoods, the search itself uses the current visit points, and only the neighbor it selects is evaluated with the BCD. The first four neighborhoods described below have no tunable parameters, while the parameters of the remaining neighborhoods are made explicit in their respective descriptions.
\begin{description}
\item[2-Opt and 3-Opt neighborhoods:] In a $k$-Opt neighborhood~\cite{Lin1965}, neighbors are generated by removing $k$ edges from the tour and reconnecting the resulting $k$ consecutive segments of cities in all possible ways to form feasible tours. For 2-Opt, $k=2$ and for 3-Opt, $k=3$. To avoid trivial or degenerate moves, edges that are consecutive along the tour are excluded from being removed. When generating neighbors for a higher-order neighborhood (e.g., 3-Opt), moves that have already been considered in lower-order neighborhoods (e.g., 2-Opt) are also excluded to avoid redundancy. The number of neighbors is $O(m^2)$ for 2-Opt and $O(m^3)$ for 3-Opt.

\item[Double-Bridge neighborhood:] The 4-Opt neighborhood~\cite{aarts2003}, as we have described above, generates neighbors by removing four edges from a tour and reconnecting the resulting segments in all possible ways; however, exhaustive enumeration is computationally expensive ($O(m^4)$). The double-bridge move is a specific 4-Opt move that partitions the tour into four consecutive segments and applies a fixed non-sequential reconnection pattern: the first segment remains unchanged, the second and fourth segments are swapped, and the third segment is preserved in its original position. This structured perturbation is widely used as a diversification mechanism, as it enables efficient escape from local optima. In our implementation, the double-bridge neighborhood is not explored by explicitly generating all possible moves. Instead, we adopt the approach of Glover~\cite{Glover1996}, which implicitly searches over a restricted but expressive class of 4-Opt moves, including double-bridge configurations and crossing variants. The search is modeled as a shortest-path problem in a directed acyclic graph, where nodes represent pairs of edges of the current tour. The best move is obtained by computing a shortest path from a source to a target node in this graph, which yields the optimal double-bridge move in $O(m^2)$ time. The arc costs of this graph are the distances between the current visit points, so the move found is optimal for these points; the BCD is then applied to the resulting permutation, which replaces the current solution if it is shorter. Importantly, no additional parameters are required in our implementation. Although Glover~\cite{Glover1996} also proposes an accelerated variant with $O(m)$ complexity based on candidate lists, this would introduce a tunable parameter controlling the maximum number of candidate edges considered. Nonetheless, this acceleration was not implemented in our case.

\item[Swap neighborhood:] Two non-adjacent cities are swapped in the current permutation. The number of neighbors is $O(m^2)$, as it considers all pairs of non-adjacent cities.

\item[Insertion neighborhood:] A city is relocated to a different position in the current permutation. The number of neighbors is also $O(m^2)$, since each city can be relocated to any of the other positions in the tour.

\item[String Relocation:] A sequence of $2 \le L \le m-1$ consecutive cities is relocated in the current permutation. The parameter $L$ and the option to insert the reversed sequence are both subject to calibration. For a fixed sequence length $L$, the neighborhood contains exactly $(m-L+1)(m-L)=O(m^2)$ neighbors. In our implementation the neighborhood relocates segments of two to $L$ consecutive cities, so every length in this range is considered at each call.

\item[Balas-Simonetti:] This neighborhood~\cite{Balas1999,Balas2001} depends on a parameter $2 \le k \le m$. A permutation of $i = (i_1, \dots, i_m)$ belongs to the neighborhood if, in the permutation, each city $i_j$ precedes all cities $i_\ell$ with $\ell \ge j + k$. Observe that the neighborhood defined by parameter $k$ is included in the neighborhood defined by $k+1$. Explicitly generating all possible neighbors of this type is computationally impractical, as the number of permutations satisfying the $k$-precedence constraints grows combinatorially with $m$ and $k$. Instead, the best neighbor for the current visit points can be efficiently identified by solving a shortest-path problem on a specific acyclic auxiliary graph $G^*$~\cite{Simonetti1996}, whose arc lengths are the distances between these points; the BCD is then applied only to the permutation given by the shortest path. The graph has $m+1$ layers: the first $m$ correspond to positions in the tour, and the $(m+1)$-th layer serves as a terminal layer representing the end of the tour. Each layer contains exactly $(k+1) 2^{k-2}$ nodes, and the number of arcs between consecutive layers is bounded by $k(k+1) 2^{k-2}$. The total number of arcs in the graph is therefore at most $m \cdot k(k+1) 2^{k-2}$. Since $G^*$ is acyclic, the shortest-path problem can be solved in time proportional to the number of nodes plus arcs, i.e., $O(m \cdot k^2 \cdot 2^k)$. In practice, if $k \ll m$ and considered constant, the method scales linearly with $m$.

\item[Gutin neighborhood:] This neighborhood~\cite{Gutin2005,Gutin2007} depends on a set $Z = \{z_1,\dots,z_{|Z|}\} \subseteq \{1,\dots,m\}$ of pairwise non-consecutive indices, meaning that $z_{i+1} > z_i + 1$ for all $i = 1,\dots,|Z|-1$ and that $1$ and $m$ cannot both belong to $Z$. It consists of the permutations obtained by removing the cities $\{ i_j : j \in Z \}$ and reassigning them to the resulting empty positions. This neighborhood contains $|Z|!$ possible neighbors. However, if the visit points of the non-removed cities remain fixed, the best neighbor can be found by solving an assignment problem on a $|Z| \times |Z|$ cost matrix $C = (c_{ab})$. Each entry $c_{ab}$ represents the cost of inserting the removed city $i_{z_a}$ into the empty position $z_b$. The conditions imposed on the set $Z$ ensure that the cities in the positions immediately before and after $z_b$ were not removed; that is, they are known. Therefore, computing $c_{ab}$ corresponds to solving problem~\eqref{bcdsubprob3} in which \(p\) and \(q\) correspond to the fixed visit points of these cities. The cost of solving this problem, as described in Section~\ref{sec3}, depends on the number of vertices of polygon~\(P_{z_a}\). Since this cost does not depend on the number of cities, we consider it an \(O(1)\) cost and the whole matrix \(C\) can be computed in $O(|Z|^2)$ time. The assignment problem can be solved exactly in $O(|Z|^3)$ time. Observe that, once the assignment problem has been solved, the BCD method is used to compute the visit points for the obtained permutation. The cardinality of $Z$ is a parameter that can be calibrated, while the set $Z$ itself is drawn at random among the admissible sets of indices at each call of the neighborhood.
\end{description}

A single BCD call is inexpensive, but in the neighborhoods explored exhaustively the number of calls grows with the size of the neighborhood, so we adopt two strategies that keep the local search efficient on the largest instances. First, the BCD is truncated after three cycles during the local search. Once the VND stops, the visit points of the resulting solution are recomputed with the BCD with its limit raised from three to 100 cycles, which in practice lets it run until its stopping criterion is met. Second, the 2-Opt, 3-Opt, swap, insertion and string relocation neighborhoods are explored with candidate lists~\cite{Helsgaun2000,Lin1973}: for each city we keep its $K$ nearest polygons, and only the moves that insert an edge between a city and one of its candidates are considered. The number of neighbors examined then drops from $O(m^2)$ to $O(mK)$ for the 2-Opt, swap, insertion and string relocation neighborhoods, and from $O(m^3)$ to $O(mK^2)$ for the 3-Opt neighborhood, with the list size $K$ left to calibration.

During the ILS perturbation step we sequentially apply a predefined number of double-bridge moves, with the edges for each move selected at random. We choose the double-bridge move over any of the other ones because of its ability to escape local minima. After having applied this sequence of moves, we recompute the visit points using the BCD method. The total number of double-bridge moves that constitute the perturbation is a parameter subject to calibration. If the local search returns the same solution in three consecutive iterations, we take this as a sign that the double-bridge perturbation is not proving effective. In that case, the next iteration applies a different perturbation: the NN heuristic is applied from a randomly chosen city and the visit points are recomputed with the BCD, which restarts the search from a different region of the solution space. The ILS stops after a prescribed number of consecutive iterations without improving the best solution found, or when a CPU time limit is reached.

\section{Numerical experiments} \label{sec5}

In this section, we describe the numerical experiments. To the best of our knowledge, there are no benchmark instance sets in the literature for the TSPN with neighborhoods defined by non-convex polygons. Therefore, in Section~\ref{sec51}, we describe how GTSP instances were adapted to construct TSPN instances with neighborhoods given by non-convex polygons. In Section~\ref{sec52}, we analyze the practical usefulness of the additional constraints~(\ref{gac1})--(\ref{gac6}). In Section~\ref{sec53}, we evaluate the robustness and efficiency of the BCD method for computing the visit points of a given sequence. In Section~\ref{sec54}, we describe the use of Irace to tune the parameters of the ILS-BCD. In Section~\ref{sec55}, we present the performance of the ILS-BCD. Finally, in Section~\ref{sec56}, we examine substantially larger instances and the limits of the method.

The BCD method described in Section~\ref{sec3}, as well as the ILS-BCD method and the neighborhoods described in Section~\ref{sec4}, were implemented in Fortran. The codes were compiled using the GNU Fortran compiler GCC (version 12.3.0) with the flags -O3 -march=native -ffast-math -flto=auto, except for the standalone BCD of Section~\ref{sec53}, which was compiled with -O3. The model~(\ref{g1})--(\ref{g24}) and its variants were solved using IBM ILOG CPLEX Optimization Studio 22.1.2 with default parameter settings. All experiments were conducted on a computer equipped with an Intel Core i9-13900 processor (5.6 GHz) and 128 GB of DDR4 RAM at 2400 MHz, running Ubuntu 23.04. The instances, the source code of the methods and the results reported in this section are available at \url{https://github.com/pcrertel/ber2026-tspn}.

\subsection{Generating TSPN instances with non-convex polygonal neighborhoods} \label{sec51}

\begin{figure}[ht!]
\begin{center}
\begin{tabular}{cc}
\includegraphics[scale=0.9]{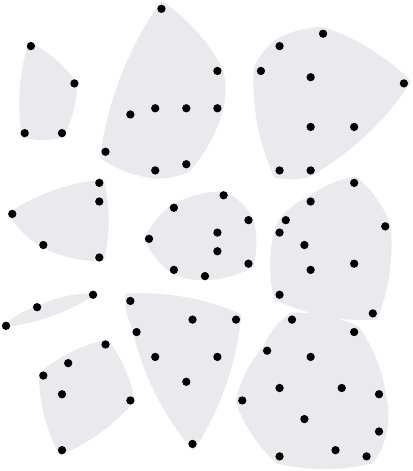} &
\includegraphics[scale=0.9]{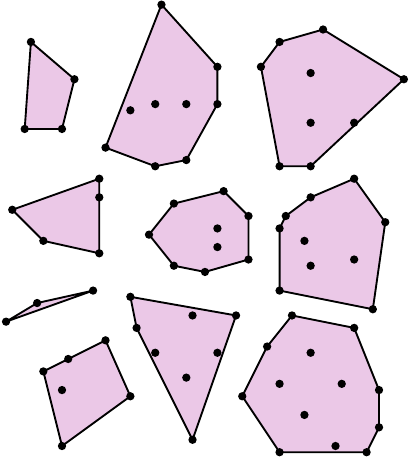} \\ [2mm]
(a) GTSP instance & (b) Convex hull of the GTSP clusters \\ [2mm]
\includegraphics[scale=0.9]{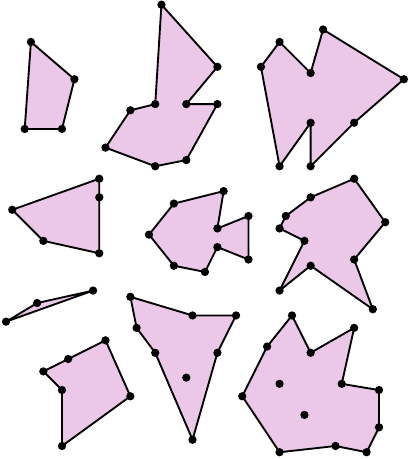} &
\includegraphics[scale=0.9]{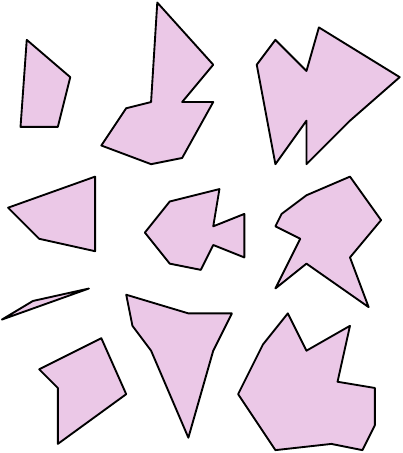} \\ [2mm]
(c) Convex hull reduction with $\hat k = 3$ & (d) Resulting TSPN instance
\end{tabular}
\end{center}
\caption{Example of the conversion of a GTSP instance from the MOM Lib~\cite{Mestria2013} library to a TSPN instance.}
\label{fig3}
\end{figure}

We construct TSPN instances from Generalized Traveling Salesman Problem (GTSP) instances. In the GTSP, each city is represented by a finite set of points; see Figure~\ref{fig3}(a) for an example containing ten cities. To convert each GTSP instance into a TSPN one, we transform each set of points---representing a city in the GTSP---into a neighborhood. Sets consisting of a single point remain unchanged. Sets formed by two or more collinear points are converted into a line segment whose endpoints are given by the pair of points that are farthest apart from one another. For sets with three or more non-collinear points, we first compute the convex hull, obtaining a polygon that contains all points in the set (see Figure~\ref{fig3}(b)). We then iteratively modify the convex hull to produce a non-convex polygon, where each modification consists of removing a triangular slice from the current polygon. Let $\hat k$ denote an upper bound on the number of modifications allowed (we set $\hat k = 3$ in our implementation). Let $k$ be the minimum between $\hat k$ and the number of interior points, i.e., points in the set that are not vertices of the convex hull. Each modification proceeds as follows. First, we identify the largest edge (say, $\overline{ab}$) of the current polygon. We then examine the interior points and form a set of candidates by selecting those interior points $p$ for which the angle between the vector $p-g$ and the horizontal axis lies between the corresponding angles of vertices $a$ and $b$, where $g$ denotes the geometric center of the GTSP cluster associated with the current city. To efficiently identify such candidates, we rely on a pre-sorted angular list. If the candidate set is nonempty, we select the candidate point $c$ with minimum Euclidean distance to $\overline{ab}$, remove the edge $\overline{ab}$, and replace it with the two edges $\overline{ac}$ and $\overline{cb}$. Otherwise, we proceed to the next largest edge of the polygon. We repeat this vertex-insertion and polygon-update procedure a total of $k$ times. This process yields a simple polygon, in general non-convex, containing all points associated with the corresponding city in the original GTSP instance. The final TSPN instance, generated from the GTSP instance in Figure~\ref{fig3}(a), is depicted in Figure~\ref{fig3}(d). 

Using this procedure, we converted 235 GTSP instances into TSPN instances with polygonal neighborhoods. Of these, 45 were taken from GTSP Lib~\cite{Fischetti1997} and 190 from MOM Lib~\cite{Mestria2013}. The selected instances comprise all instances available in these libraries with up to 100 cities for which the cities' geographic coordinates are provided, as these coordinates are required for the conversion. The resulting instances contain between 8 and 100 polygons, with the total number of vertices ranging from 30 to 1{,}163. For the instances considered, model~(\ref{g1})--(\ref{gac6}) has up to approximately 10{,}000 binary variables, 40{,}000 continuous variables, 100{,}000 linear constraints, and 10{,}000 quadratic constraints. To represent an instance using model~(\ref{g1})--(\ref{gac6}), each nonconvex polygonal neighborhood must be represented as a union of convex polygons. To this end, we opted to partition each neighborhood into convex polygons. This partitioning was performed using the \texttt{optimal\_convex\_partition\_2} routine from the CGAL library~\cite{cgal:h-pp2-25b}. After partitioning the nonconvex neighborhoods, the total number of convex polygons ranges from 12 to 297, while the total number of vertices of these convex subpolygons ranges from 34 to 1{,}531. Tables~S1 and~S2 in the supplementary material provide detailed information on the instances adapted from GTSP Lib and MOM Lib, respectively.

\subsection{Evaluation of the additional constraints based on norm inequalities} \label{sec52}

To assess whether the additional constraints~(\ref{gac1})--(\ref{gac6}) improve the efficiency of CPLEX in solving the MIQCP model~(\ref{g1})--(\ref{g24}), we solved the 87 TSPN instances with at most 15 polygons using both formulations, namely~(\ref{g1})--(\ref{g24}) and~(\ref{g1})--(\ref{gac6}). All CPLEX parameters were kept at their default values, and a CPU time limit of one hour was imposed for each run. We constructed an initial solution using the nearest neighbor (NN) constructive heuristic. To this end, we represent each city by the centroid of its neighborhood. The method starts the tour from an arbitrary city and, at each iteration, appends the city closest to the current last city among those not yet visited. Since the tour produced depends on the choice of the starting city, we run the heuristic $m$ times, each time starting from a different one of the $m$ cities. Among the $m$ tours obtained, we select the one with the shortest total length. The solution $(s, w)$ provided to CPLEX as an initial solution is defined as follows: $s$ is the sequence constructed by the NN heuristic described above, and $w$ corresponds to the centroids.

For the baseline model~(\ref{g1})--(\ref{g24}), CPLEX proved optimality for 80 of the 87 instances. When applied to the model~(\ref{g1})--(\ref{gac6}), CPLEX proved optimality for 78 of them, all of which belong to the set of 80 instances solved by the baseline model; the two exceptions are 12nrw1379-3x4 and 15pr76-3x5, for which the baseline model proved optimality in 3{,}481.37s and 1{,}657.77s, respectively, whereas the model with the additional constraints reached the time limit with gaps of 18.07\% and 15.35\%. Table~S3 in the supplementary material reports, for each model, the tour length of the solution returned by CPLEX (when optimality is proven) or the interval defined by the lower and upper bounds along with the corresponding gap. The table also presents the number of iterations, the number of nodes explored in the branch-and-bound tree, and the CPU time in seconds.

For comparison purposes, the 87 instances were partitioned into two subsets: (a) the set of 78 instances for which CPLEX proved optimality for both models~(\ref{g1})--(\ref{g24}) and~(\ref{g1})--(\ref{gac6}); and (b) the remaining 9 instances, for which CPLEX did not prove optimality for at least one of the two models. For subset (a), where both models return a certified optimal tour and therefore the same objective function value, the comparison is made in terms of CPU time. To this end, we applied the Wilcoxon test~\cite{Wilcoxon1945}, with a significance level of $\alpha = 0.05$, to test the null hypothesis that ``the samples from the two models come from the same distribution'' or, equivalently, that ``the differences between paired samples follow a symmetric distribution around zero.'' Subset (b), with only 9 instances, is too small to support a meaningful statistical test; it is therefore analyzed instance by instance.

Table~\ref{tab1} summarizes the computational comparison between the two formulations. For the 78 instances solved to proven optimality by both models, adding constraints~(\ref{gac1})--(\ref{gac6}) reduces mean CPU time, median CPU time, and the mean number of branch-and-bound nodes by $21.82\%$, $31.28\%$, and $37.73\%$, respectively. The formulation with additional constraints is faster on 48 of the 78 instances, with four ties; excluding these ties, the Wilcoxon signed-rank test on the remaining 74 paired observations yields a $p$-value of $1.83\times10^{-3}$, providing statistical evidence of improved computational efficiency at the $\alpha=0.05$ significance level. The results are less conclusive for the seven instances on which neither formulation proves optimality: the baseline model achieves a smaller final gap in five cases, owing to tighter lower bounds, whereas the formulation with additional constraints finds a shorter incumbent tour in two cases and an equally short tour in the other five. On the two remaining instances, only the baseline formulation proves optimality. For 15pr76-3x5, the incumbent returned by the formulation with additional constraints is marginally shorter than the tour certified as optimal by the baseline model, a difference within CPLEX's default relative-gap tolerance; for 12nrw1379-3x4, it is $0.02\%$ longer. Overall, the additional constraints improve the computational effort required to prove optimality on instances solved by both formulations, but their impact on instances that remain unsolved within the time limit is mixed: they can yield slightly better incumbent solutions, but do not consistently improve the final optimality gap.

\begin{table}[ht!]
\centering
\small
\setlength{\tabcolsep}{5pt}
\renewcommand{\arraystretch}{1.12}
\begin{tabular}{@{}lrrr@{}}
\toprule
& MIQCP~(\ref{g1})--(\ref{g24}) & MIQCP~(\ref{g1})--(\ref{gac6}) & Reduction (\%) \\
\midrule
\multicolumn{4}{@{}l}{\textit{(a) Both models prove optimality (78 instances)}} \\
\addlinespace[2pt]
Mean CPU time (s)       & 177.09      & 138.45      & 21.82 \\
Median CPU time (s)     & 35.29       & 24.25       & 31.28 \\
Mean B\&B nodes         & 1{,}886{,}146 & 1{,}174{,}511 & 37.73 \\
Instances solved faster & 26          & 48          & --- \\
\midrule
\multicolumn{4}{@{}l}{\textit{(b) Neither model proves optimality (7 instances)}} \\
\addlinespace[2pt]
Mean final gap (\%)          & 39.05 & 41.30 & --- \\
Instances with smaller gap   & 5     & 2     & --- \\
Instances with shorter tour  & 0     & 2     & --- \\
\bottomrule
\end{tabular}
\caption{Summary of the comparison between the baseline MIQCP formulation~(\ref{g1})--(\ref{g24}) and the formulation augmented with constraints~(\ref{gac1})--(\ref{gac6}), based on the 87 instances with at most 15 polygons (see Table~S3 in the supplementary material). Subset (a) comprises the 78 instances for which both models prove optimality; subset (b) comprises the seven instances for which neither model proves optimality. The four CPU-time ties in subset (a) are excluded from the counts of instances solved faster. The two remaining instances, on which only the baseline model proves optimality, are discussed in the text. Reductions in CPU time and B\&B nodes are computed from the corresponding means.}
\label{tab1}
\end{table}

\subsection{Analysis of the use of BCD for computing visit points of a given sequence} \label{sec53}

Given a sequence $s$, we propose computing the visit points $w$ using the Block Coordinate Descent (BCD) method. This problem can also be modeled as an MIQCP and solved with CPLEX. Specifically, the MIQCP for computing visit points for a given sequence corresponds exactly to model~(\ref{g1})--(\ref{gac6}), in which the variables $x_{ij}$, for $i,j=1,\dots,m$, $i \neq j$, and the variables $u_i$, for $i=2,\dots,m$, are fixed to represent the given sequence $s=(s_1,\dots,s_m)$. In particular, $x_{ij}=1$ for $(i,j) \in \{ (s_1,s_2), (s_2,s_3), \dots, (s_{m-1},s_m), (s_m,s_1) \}$, and $x_{ij}=0$ for all other pairs $(i,j)$ with $i \neq j$. The variables $u_i$ are set equal to the displacement of city $i$ in the sequence with respect to the position of city~1; if the displacement is negative, $m$ is added. We refer to this formulation as the restricted model~(\ref{g1})--(\ref{gac6}). To evaluate the performance and the quality of the solutions (visit points) produced by the BCD method, we conducted computational experiments in which, for a given sequence $s$, we solved the problem of computing the visit points using both BCD and the restricted model~(\ref{g1})--(\ref{gac6}). The ordering $s$ was constructed using the NN heuristic, as described in Section~\ref{sec52}.

For a given sequence $s$, the BCD method computes the visit points $w$ by iteratively optimizing one point at a time while keeping the others fixed. This process is repeated cyclically, meaning that each cycle $r = 1,2,\dots$ consists of $m$ consecutive iterations $t = (r-1)m+1, (r-1)m+2, \dots, (r-1)m+m$. At each iteration $t$, the visit point $p_{i_{\ell}}$, with $\ell = \text{mod}(t,m)$, is computed by solving subproblem~\eqref{bcdsubprob}. Let $w^r$ denote the visit points computed at cycle~$r$, and let $f(s,w^r)$ be the length of the tour defined by the sequence $s$ and the visit points $w^r$. We terminate the BCD iterative process at cycle $r$ when the tour length associated with $w^r$ does not exhibit a significant improvement relative to that of the previous cycle, i.e.,
\begin{equation*} \label{stopcrit}
f(s,w^r) \not\leq f(s,w^{r-1}) - f_{\mathrm{tol}} \, \max\{1, |f(s,w^{r-1})|\},
\end{equation*}
where $f_{\mathrm{tol}}=\sqrt{\epsilon_{\mathrm{mach}}}$ and $\epsilon_{\mathrm{mach}} = 10^{-16}$. For $r = 0$, $w^r$ corresponds to the centroids; thus, both BCD and CPLEX start from the same initial guess for the visit points.

All CPLEX parameters were kept at their default values, and a CPU time limit of one hour was imposed for each run. Tables~S4 and~S5 in the supplementary material report the instance-level results for the 45 instances adapted from the GTSP Lib and the 190 instances adapted from the MOM Lib, respectively. Table~\ref{tab2} summarizes the results over all 235 instances, distinguishing between those for which CPLEX proved optimality and those for which it reached the time limit. CPLEX proved optimality in 215 instances. Among these, BCD obtained solutions with the same tour length, up to a relative difference of $0.01\%$, in 209 instances. In the remaining six instances, BCD produced tours that were, on average, only $0.08\%$ longer than the optimal tours, with the largest difference being $0.24\%$. For these 215 instances, the average CPU time of CPLEX was 204.12 s, whereas BCD required, on average, $4.11\times10^{-5}$ seconds, i.e., more than six orders of magnitude faster. For the remaining 20 instances, CPLEX reached the one-hour time limit, with an average optimality gap of $2.50\%$. In this group, BCD found shorter tours in 17 instances, a tour of equal length in two instances, and a longer tour in only one instance. The shorter tours were, on average, $1.21\%$ shorter than the corresponding CPLEX solutions, while the worse solution was only 0.07\% longer. BCD required an average of $2.58\times10^{-4}$ s per instance, whereas CPLEX did not finish within the one-hour limit. These results show that BCD produces solutions of comparable quality at a negligible computational cost. This makes it particularly suitable for use within a heuristic framework, where the rapid evaluation of many candidate sequences is essential.

\begin{table}[ht!]
\centering
\small
\renewcommand{\arraystretch}{1.15}
\setlength{\tabcolsep}{4pt}
\resizebox{\textwidth}{!}{%
\begin{tabular}{lrrrrrrr}
\toprule
\multirow{2}{*}{\raisebox{-0.10ex}{Instance}}
& \multirow{2}{*}{\raisebox{-0.10ex}{\# Inst.}}
& \multicolumn{3}{c}{CPLEX}
& \multicolumn{2}{c}{BCD}
& \multirow{2}{*}{\raisebox{-0.10ex}{Rel. diff. (\%)}} \\
\cmidrule(lr){3-5} \cmidrule(lr){6-7}
& & \multicolumn{1}{c}{\# iterations}
& \multicolumn{1}{c}{\# B\&B nodes}
& \multicolumn{1}{c}{CPU (s)}
& \multicolumn{1}{c}{\# cycles}
& \multicolumn{1}{c}{CPU (s)}
& \\
\midrule
Optimality proved
& 215 & 29{,}245 & 1{,}056 & 204.12 & 18.6
& $4.11\times10^{-5}$ & 0.002 \\
Time limit
& 20 & 797{,}020 & 30{,}135 & 3{,}608.03 & 27.1
& $2.58\times10^{-4}$ & -1.028 \\
\midrule
All
& 235 & 94{,}587 & 3{,}531 & 493.82 & 19.3
& $5.95\times10^{-5}$ & -0.088 \\
\bottomrule
\end{tabular}%
}
\caption{Summary of the computational results for computing the visit points of a fixed sequence using CPLEX and BCD. Each entry reports the mean over the instances in the corresponding group. ``Optimality proved'' denotes instances for which CPLEX proved optimality, whereas ``Time limit'' denotes instances for which CPLEX reached the one-hour time limit. The relative difference is computed as $100\% \times (f_{\mathrm{BCD}}-f_{\mathrm{CPLEX}})/f_{\mathrm{CPLEX}}$. Negative values indicate shorter tours obtained by BCD. Instances with zero CPLEX tour length are excluded from this calculation.}
\label{tab2}
\end{table}

The two approaches also differ in the feasibility of the visit points they return. BCD produces visit points that belong to their respective neighborhoods up to machine precision, whereas CPLEX may return points that lie slightly outside the neighborhoods, within the solver's feasibility tolerance. This explains the small nonzero differences in tour length reported in Tables~S4 and~S5, which often favor CPLEX even when both methods produce solutions of comparable quality. A point lying marginally outside its neighborhood can slightly shorten the tour, and such deviations may accumulate over the tour. Consequently, these small differences should be interpreted in light of CPLEX's feasibility tolerance rather than as evidence of a meaningful difference in solution quality.

The comparison between computing the visit points using BCD and solving the restricted MIQCP model with an exact method is intended only to demonstrate that BCD enables the development of efficient heuristic methods for the TSPN due to its low computational cost and high-quality solutions. On the instances considered, the BCD method was able to find optimal visit points in most cases. However, it is not guaranteed to do so in general, as illustrated in Figure~\ref{fig2}: BCD converges to a coordinate-wise minimum, which need not be a global minimum of the full problem. Section~\ref{sec531} isolates the two hypotheses under which the guarantee does hold and measures how much they matter in practice. Furthermore, since the use of BCD to compute visit points for a given sequence is embedded within a heuristic framework, we adopted its fastest variant, as described above. This choice allows the heuristic, within a fixed time limit, to evaluate a larger number of sequences, which is a more critical factor for the overall performance of the method. For the sake of completeness, another variant that we tested in the experiments is described in Section~\ref{sec532}.

\subsubsection{Effect of non-convexity and overlap of the neighborhoods on BCD} \label{sec531}

The two hypotheses stated in Section~\ref{sec3} under which BCD is guaranteed to find optimal visit points for a fixed sequence are that all neighborhoods are convex and that consecutive neighborhoods in the sequence are disjoint. These conditions do not hold simultaneously for all instances in the original test set, and, as discussed in Section~\ref{sec53}, BCD may therefore fail to find the optimal visit points in seven of them. To investigate the effect of non-convexity of the neighborhoods and of the overlap between them, we constructed a modified set of TSPN instances from the 235 GTSP instances using the procedure described in Section~\ref{sec51}, stopping after replacing each neighborhood with the convex hull of its corresponding GTSP cluster (Figure~\ref{fig3}(b)). Thus, every neighborhood in the modified instances is a single convex polygon. For each instance, we used the sequence generated by the nearest neighbor (NN) heuristic and computed the visit points using both BCD and CPLEX applied to the restricted model~(\ref{g1})--(\ref{gac6}). Both methods were initialized at the centroids of the neighborhoods. Since all neighborhoods are convex, the restricted problem reduces to a second-order cone program (SOCP), and CPLEX proved optimality for all 235 instances.

As summarized in Table~\ref{tab3}, BCD matched the tour length obtained by CPLEX in all 235 instances, within a relative tolerance of $0.01\%$, and does so in $5.08 \times 10^{-5}$ seconds on average against $35.65$ seconds for CPLEX, nearly six orders of magnitude faster and essentially the same time it takes on the non-convex instances. In 182 of the 235 instances, every pair of consecutive neighborhoods of the sequence is disjoint, so both hypotheses hold and every coordinate-wise minimum is a global minimizer. BCD needs $13.84$ cycles on average, with a median of $10.5$, a minimum of 3 and a maximum of 67, and in neither solution does any pair of consecutive visit points coincide. In the remaining 53 instances at least one pair of consecutive neighborhoods of the sequence intersects and differentiability is no longer guaranteed. In 18 of these instances the BCD solution has at least one pair of coincident consecutive visit points, and the CPLEX solution has one in 19. BCD still reaches the tour length of CPLEX in every one of them, at the cost of more cycles, $38.64$ on average. This average is dominated by six instances whose neighborhoods share a common point: there every consecutive pair of the sequence overlaps, all visit points collapse onto that point, the optimal tour has length zero, and BCD needs $171.67$ cycles on average, up to $340$. Over the other 47 instances it needs $21.66$ cycles on average, with a median of 17 and a maximum of 71, against $13.84$ in the disjoint group. Tables~S6 and~S7 in the supplementary material report the instance-level results, while Table~\ref{tab3} summarizes the results by group.

\begin{table}[ht!]
\centering
\small
\renewcommand{\arraystretch}{1.15}
\setlength{\tabcolsep}{4pt}
\begin{tabular}{lrrrrrr}
\toprule
\multirow{2}{*}{\raisebox{-0.10ex}{Instance}}
& \multirow{2}{*}{\raisebox{-0.10ex}{\# Inst.}}
& \multicolumn{2}{c}{CPLEX}
& \multicolumn{2}{c}{BCD}
& \multirow{2}{*}{\raisebox{-0.10ex}{Rel. diff. (\%)}} \\
\cmidrule(lr){3-4} \cmidrule(lr){5-6}
& & \multicolumn{1}{c}{\# iterations}
& \multicolumn{1}{c}{CPU (s)}
& \multicolumn{1}{c}{\# cycles}
& \multicolumn{1}{c}{CPU (s)}
& \\
\midrule
Disjoint                      & 182 & 2{,}791 & 39.36 &  13.84 & $4.2\times10^{-5}$ & 0.00 \\
Intersecting (positive length) &  47 & 3{,}971 & 25.85 &  21.66 & $8.5\times10^{-5}$ & 0.00 \\
Intersecting (zero length)     &   6 &       0 &  0.06 & 171.67 & $6.7\times10^{-5}$ & --- \\
\midrule
All                           & 235 & 2{,}956 & 35.65 &  19.43 & $5.1\times10^{-5}$ & 0.00 \\
\bottomrule
\end{tabular}
\caption{Summary of the computational results for computing the visit points of a fixed sequence using CPLEX and BCD on the 235 instances obtained by replacing each GTSP cluster with its convex hull. The instances are grouped according to whether consecutive neighborhoods in the sequence are disjoint or intersecting. ``Intersecting (positive length)'' denotes instances with at least one intersecting consecutive pair and a positive optimal tour length. ``Intersecting (zero length)'' denotes instances in which all neighborhoods share a common point, yielding a zero-length optimal tour. Each entry reports the mean over the instances in the corresponding group. The relative difference is computed as $100\% \times (f_{\mathrm{BCD}}-f_{\mathrm{CPLEX}})/f_{\mathrm{CPLEX}}$. Negative values indicate shorter tours obtained by BCD. For the six zero-length instances, the relative difference is undefined. For all other instances, the relative difference rounds to $0.00\%$.} 
\label{tab3}
\end{table}

\subsubsection{Multi-start BCD method} \label{sec532}

In the experiment described in Section~\ref{sec53}, seven out of 235 instances yielded tours from BCD that were slightly worse than those obtained by CPLEX applied to the restricted model~(\ref{g1})--(\ref{gac6}), an instance of the convergence to a non-global coordinate-wise minimum discussed in Section~\ref{sec3}. We apply here the multi-start strategy mentioned there: BCD is run $n_{\mathrm{trial}}$ times and the shortest of the $n_{\mathrm{trial}}$ tours obtained is kept. The first run starts from the centroids, exactly as in the main experiment. In each of the remaining $n_{\mathrm{trial}}-1$ runs, every visit point is drawn uniformly at random on a randomly chosen edge of its neighborhood, which guarantees a feasible starting point.

Increasing $n_{\mathrm{trial}}$ progressively reduces the number of instances in which BCD returns a longer tour than CPLEX: three of the seven instances are brought to the quality of CPLEX with $n_{\mathrm{trial}}=2$, five with $n_{\mathrm{trial}}=3$, six with $n_{\mathrm{trial}}=5$ and $10$, and all seven with $n_{\mathrm{trial}}=20$, so that BCD no longer returns a longer tour in any of the 235 instances. The price is the CPU time, which grows roughly in proportion to $n_{\mathrm{trial}}$: from $5.95 \times 10^{-5}$ seconds on average with $n_{\mathrm{trial}}=1$ to $1.04 \times 10^{-3}$ seconds with $n_{\mathrm{trial}}=20$. Within a metaheuristic that computes visit points for fixed visiting sequences a huge number of times, such an increase is not justified by a gain in tour length that, in fact, never exceeds $0.24\%$ (over the seven instances whose tours were improved with the multi-start strategy). Moreover, within the ILS, BCD is not started from the centroids of the neighborhoods, but from the visit points of the current solution of the ILS or of the VND, which are expected to be already close to the optimal visit points of the sequence being evaluated; the benefit of restarting the method from other points is therefore expected to be even smaller in that context. This is why the single-start variant is the one embedded in the ILS.

\subsection{Parameter tuning of the ILS and neighborhood structures} \label{sec54}

The ILS-BCD of Section~\ref{sec4} has eleven parameters to be calibrated. Three of them govern the acceptance criterion: the initial threshold $\eta_0 \in \{0.01, 0.02, \ldots, 0.09\}$, under which a solution is accepted if its cost exceeds that of the best solution found by at most a fraction $\eta$; the factor $h \in \{0.80, 0.85, 0.90, 0.95, 0.99, 1\}$ by which the threshold is multiplied at each update; and the number $n_{\mathrm{stag}} \in \{5, 10, 25, 50, 100\}$ of consecutive iterations without improvement between two updates of the threshold. Six parameters concern the variable neighborhood descent: the order in which its eight neighborhoods are visited, either fixed or drawn at random at each call; the improvement strategy, first or best improvement, of the 2-opt, 3-opt, swap, insertion and string relocation neighborhoods; the parameter $k \in \{2, 3\}$ of the Balas--Simonetti neighborhood; the number $|Z| \in \{5, 10, \ldots, 30\}$ of tour positions reassigned by solving an assignment problem in the neighborhood of Gutin; the maximum length $L \in \{2, 3, 4\}$ of the strings moved by the string relocation, which relocates segments of two up to $L$ consecutive cities of the tour; and whether a relocated string may also be inserted in reverse. Finally, the perturbation applies $\ell \in \{1, \ldots, 5\}$ double-bridge moves in sequence, each one to the tour produced by the previous one, and the size $K \in \{5, 10, \ldots, 30\}$ of the candidate lists bounds the width of the 2-opt, 3-opt, swap, insertion and string relocation searches. Together, these sets define a space of 2{,}332{,}800 configurations.

\begin{figure}[ht!]
\begin{center}
\begin{tabular}{ccc}
\includegraphics[scale=0.6]{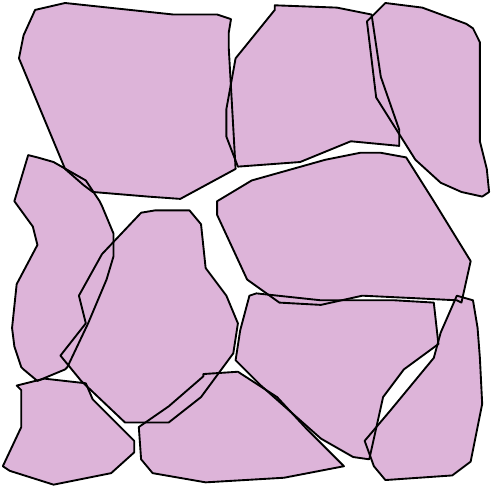} &
\includegraphics[scale=0.6]{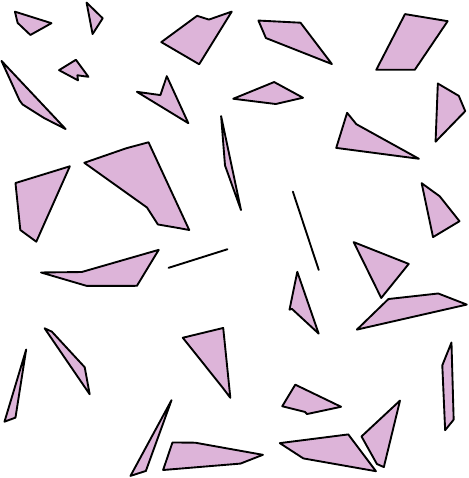} &
\includegraphics[scale=0.6]{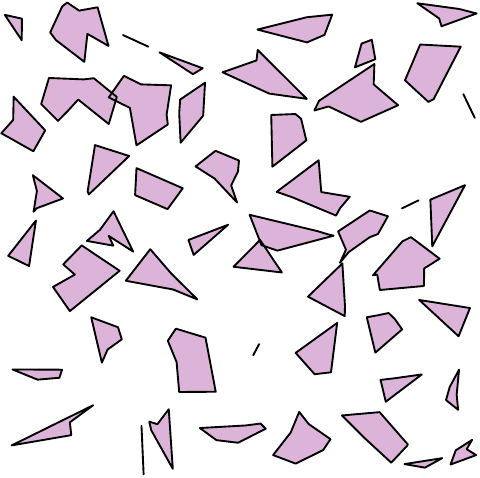} \\ [2mm]
(a) $m = 10$ & (b) $m = 30$ & (c) $m = 50$ \\ [2mm]
\includegraphics[scale=0.6]{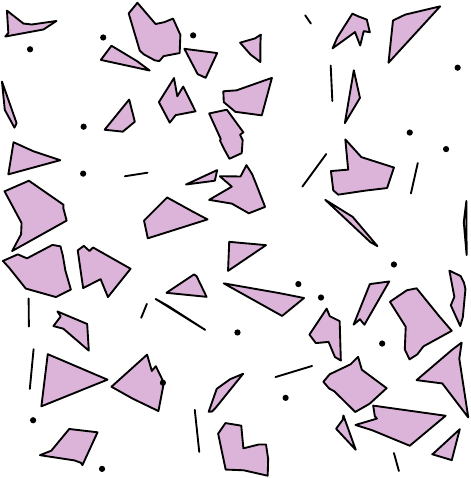} &
\includegraphics[scale=0.6]{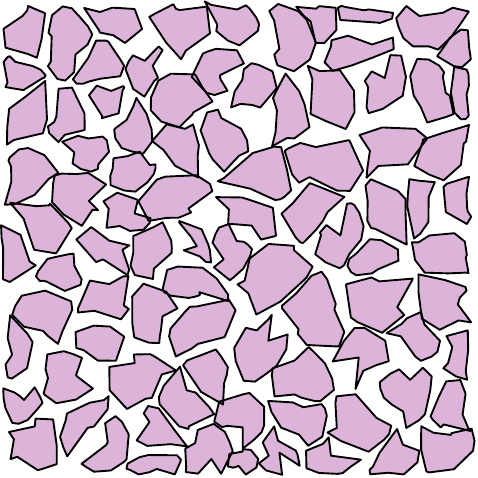} &
\includegraphics[scale=0.6]{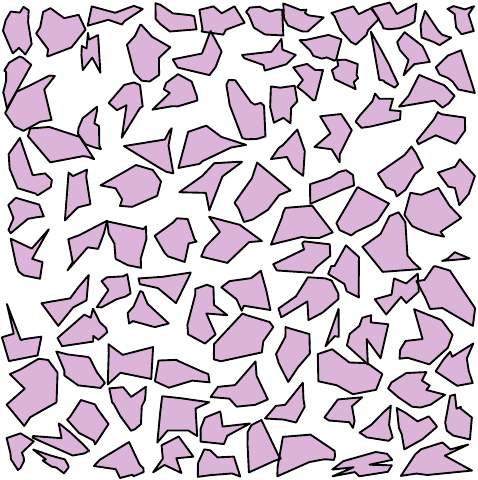} \\ [2mm]
(d) $m = 70$ & (e) $m = 90$ & (f) $m = 100$
\end{tabular}
\end{center}
\caption{Six of the 60 instances generated for the calibration, in increasing order of the number of neighborhoods $m$. Each neighborhood is a simple polygon with at most three reflex vertices.}
\label{fig4}
\end{figure}

For the calibration of the parameters, we generated a set of random instances, disjoint from the 235 instances converted in Section~\ref{sec51}, with a generator that reproduces the procedure by which the converted instances were obtained. To build an instance with $m$ neighborhoods on a square of side $D$, we first build $m$ clusters of points, each of which gives rise to one polygonal neighborhood, as follows. The centers of the clusters are drawn by best-candidate sampling~\cite{Mitchell1991}: each center is the best of ten candidates drawn uniformly at random in the square, namely the one whose distance to the centers already placed is the largest, which spreads the centers more evenly than uniform sampling. The number of points of each cluster is then taken from the GTSP instances of Section~\ref{sec51}: we pick at random one of them whose number of clusters lies between $0.8m$ and $1.3m$ and shuffle the list with the numbers of points of its clusters; if that instance has more than $m$ clusters, the first $m$ entries of the shuffled list are used, and otherwise the list is traversed cyclically until $m$ entries are obtained. Next, the points of each cluster are drawn uniformly in the Voronoi cell of its center, clipped to the square, and each point is displaced in a random direction by a random distance of at most one tenth of the distance from its center to the nearest other center. Points can thus slightly cross the boundary of their cell, which reproduces the small overlap between the neighborhoods of the converted instances. The coordinates are then rounded to integers. Finally, each cluster is turned into a polygon as in Section~\ref{sec51}. Since coordinates are integers, the side $D$ also sets how coarse the grid is relative to a cluster, and hence how often a cluster collapses to a segment or a point. Four sides are used, $D = 2 \times 10^{j}$ with $j \in \{2, 3, 4, 5\}$, one per order of magnitude, within the range of coordinates of the converted instances. Using this procedure, we generated 60 instances with $m \in \{10, 20, \ldots, 100\}$: five for each $m \in \{10, 20, 30, 40\}$, six for each $m \in \{50, 60, 70, 80\}$ and eight for each $m \in \{90, 100\}$. More instances are generated for the larger values of $m$, because those instances are the hardest and the most sensitive to the values of the parameters. Figure~\ref{fig4} shows six of the generated instances.

The eleven parameters were calibrated with Irace~\cite{irace} on the 60 random instances, split into a training set of 20, used by Irace, and a test set of 40, used only to evaluate the final configurations. The training set gives more weight to the larger instances, with one instance for each $m \in \{10, 20, 30, 40\}$, two for each $m \in \{50, 60, 70, 80\}$ and four for each $m \in \{90, 100\}$, and the domain sides taken in turn; the test set has four instances for each value of $m$, one per domain side. The budget was 5{,}000 runs of the ILS-BCD, each limited to 300 seconds and with its own seed, of which 4{,}695 were spent over five iterations, and configurations were eliminated by the Friedman test at the 95\% confidence level. The value returned to Irace is the tour length in normalized coordinates, those of each instance divided by their largest absolute value. Irace returned the five elite configurations of Table~\ref{tab4}, listed in the order of its final ranking, which were then run on the 40 test instances. They agree on the structure of the search: all of them visit the neighborhoods of the VND in fixed order, use first improvement and forbid the insertion of reversed strings, and four of the five start from $\eta_0 = 0.01$, set $k = 3$ in the Balas--Simonetti neighborhood and apply five double-bridge moves in the perturbation step. The candidate lists are much shorter than the longest ones offered, with $K = 10$ in three of the configurations and $K = 5$ in the other two. The remaining parameters, $h$, $n_{\mathrm{stag}}$, $L$ and $|Z|$, vary among the elite configurations.

\begin{table}[ht!]
\centering
\small
\renewcommand{\arraystretch}{1.15}
\begin{tabular}{c@{\hspace{15pt}}ccccccccc}
\toprule
Rank
& $\eta_0$
& $h$
& $n_{\mathrm{stag}}$
& $k$
& $|Z|$
& $L$
& Reverse
& $\ell$
& $K$ \\
\midrule
1 & 0.01 & 0.95 & 10 & 3 & 10 & 3 & no & 5 & 10 \\
2 & 0.01 & 0.95 & 25 & 3 & 15 & 2 & no & 5 &  5 \\
3 & 0.01 & 0.95 & 25 & 3 & 10 & 2 & no & 4 &  5 \\
4 & 0.01 & 0.80 & 25 & 3 & 15 & 2 & no & 5 & 10 \\
5 & 0.02 & 0.90 & 10 & 2 & 10 & 2 & no & 5 & 10 \\
\bottomrule
\end{tabular}
\caption{Elite configurations returned by Irace~\cite{irace}, in the order of its final ranking.}
\label{tab4}
\end{table}

On the 40 test instances the five elite configurations reach the same solution on 31 of them, including every instance with $m \le 70$. They differ only on instances with $m \ge 80$, and a Friedman test applied to the results of the five configurations on the test instances detects no significant difference among them ($p = 0.69$). Since the test set does not discriminate among them, we adopted the configuration ranked first by Irace. This configuration also ties for the largest number of test instances on which the best solution among the five is reached, 36 of the 40, together with those ranked third and fourth. Its values are $\eta_0 = 0.01$, $h = 0.95$, $n_{\mathrm{stag}} = 10$, fixed neighborhood order, first improvement, $k = 3$, $|Z| = 10$, $L = 3$ with reversed strings not allowed, $\ell = 5$ and $K = 10$, and this configuration is used in all the computational experiments in the next section.

\subsection{Evaluation of the performance of the ILS} \label{sec55}

In this section, we evaluate the ILS-BCD proposed in Section~\ref{sec4}, using the parameter configuration selected in Section~\ref{sec54}. We first consider the 235 TSPN instances obtained by the conversion procedure described in Section~\ref{sec51} and then examine an additional instance from the literature~\cite{Birgin2022}. Optimal solutions are certified for only 80 of the converted instances, those solved to proven optimality in Section~\ref{sec52}; we take them up separately below. For the set as a whole, we compare the ILS-BCD results with the best-known solutions of the corresponding GTSP instances, obtained with the solvers of~\cite{Helsgaun2015}. We refer to these as the GTSP tours and denote their lengths by $f_{\mathrm{GTSP}}$. Because every point in a GTSP cluster belongs to the polygon constructed from that cluster, any GTSP tour is also a feasible TSPN tour. Let $f_{\mathrm{TSPN}}$ denote the length of a tour obtained by the ILS-BCD. We define the relative gap as $100\% \times (f_{\mathrm{TSPN}}-f_{\mathrm{GTSP}})/f_{\mathrm{GTSP}}$, so a negative gap indicates a tour shorter than the corresponding GTSP tour. Each instance was solved in ten independent runs with different random seeds. A run terminates after 3{,}000 consecutive ILS-BCD iterations without improving the incumbent solution or after reaching a CPU time limit of 1{,}200 seconds. All reported CPU times include the construction of the initial solution. Table~S8 in the supplementary material reports the results for individual instances, while Table~\ref{tab5} summarizes the results by instance size.

\begin{table}[ht!]
\centering
\small
\renewcommand{\arraystretch}{1.15}
\setlength{\tabcolsep}{4pt}
\begin{tabular}{lrrrrrrr}
\toprule
\multirow{2}{*}{Instance size ($m$)}
& \multirow{2}{*}{\raisebox{-0.10ex}{\# Inst.}}
& \multicolumn{3}{c}{Gap (\%)}
& \multirow{2}{*}{\shortstack[r]{Improvement over\\init. solution (\%)}}
& \multirow{2}{*}{\shortstack[r]{Mean CPU\\time (s)}}
& \multirow{2}{*}{\raisebox{-0.10ex}{CV (\%)}} \\
\cmidrule(lr){3-5}
& & \multicolumn{1}{c}{Best}
& \multicolumn{1}{c}{Mean}
& \multicolumn{1}{c}{Worst}
& & & \\
\midrule
$m \le 10$         &  69 & -0.555 & -0.555 & -0.555 &  1.31 &    6.32 & 0.0000 \\
$11 \le m \le 20$  &  59 & -0.353 & -0.353 & -0.353 &  6.75 &   30.71 & 0.0000 \\
$21 \le m \le 50$  &  81 & -0.432 & -0.432 & -0.432 & 13.79 &  274.52 & 0.0000 \\
$51 \le m \le 100$ &  26 & -0.376 & -0.356 & -0.312 & 18.12 &  990.45 & 0.0237 \\
\midrule
All                & 235 & -0.441 & -0.438 & -0.433 &  8.99 &  213.77 & 0.0027 \\
\bottomrule
\end{tabular}
\caption{Summary of the ILS-BCD results on the 235 TSPN instances, grouped by the number of neighborhoods, $m$. The table reports the mean, over the instances in each group, of the best, mean and worst gaps across ten independent runs, the mean improvement of the best run over the initial solution, the mean CPU time per run, and the coefficient of variation (CV) of the tour lengths across the ten runs. The initial solution is the multi-start nearest neighbor construction followed by BCD described in Section~\ref{sec4}, which involves no randomness and is therefore the same in the ten runs. The six instances whose neighborhoods share a common point, and hence have a zero-length TSPN tour, are excluded from the gap, improvement and CV calculations.}
\label{tab5}
\end{table}

\begin{figure}[ht!]
\begin{center}
\begin{tabular}{cc}
\includegraphics[scale=0.9]{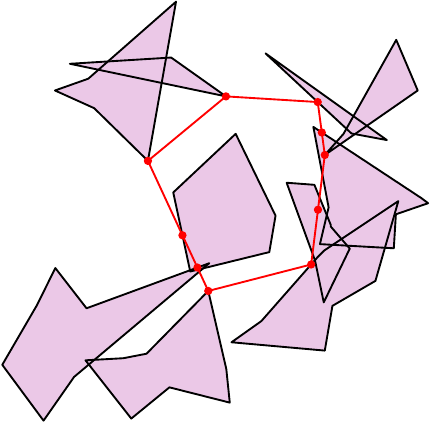} &
\includegraphics[scale=0.9]{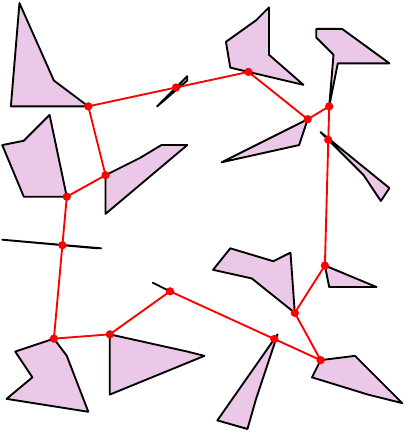} \\ [2mm]
(a) 10i60-21, gap $-2.71\%$ & (b) 16st70-4x4, gap $-1.78\%$ \\ [2mm]
\includegraphics[scale=0.9]{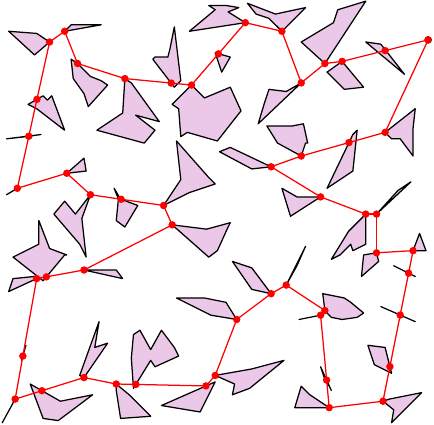} &
\includegraphics[scale=0.9]{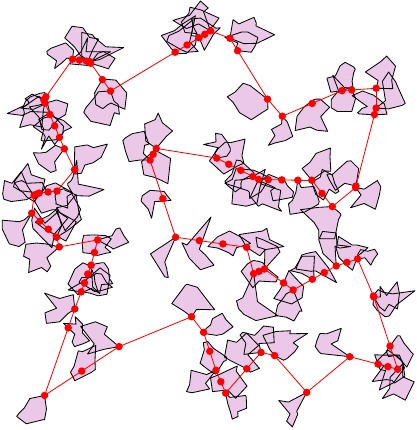} \\ [2mm]
(c) 53gil262, gap $-0.71\%$ & (d) 100i2000-604, gap $-0.42\%$
\end{tabular}
\end{center}
\caption{Best solution found by the ILS-BCD on four instances, in increasing order of the number of polygons, each with its gap to the GTSP tour. The tour is drawn in red and the visit points are the red dots.}
\label{fig5}
\end{figure}

\begin{figure}[ht!]
\centering
\includegraphics[scale=0.85]{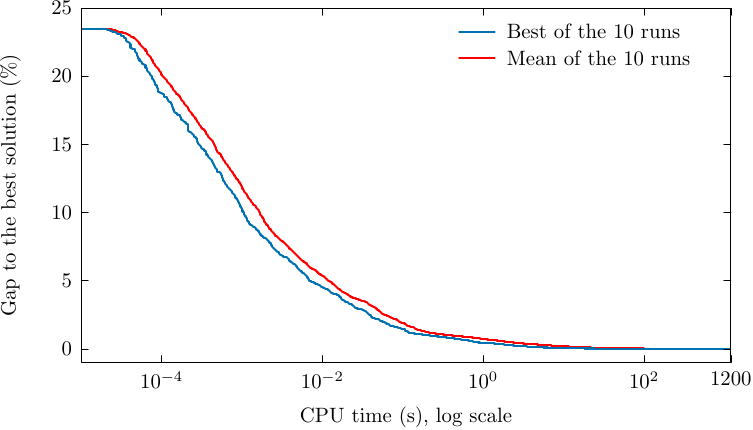}
\caption{Improvement of the incumbent solution over time: for each instant, the value plotted is the mean, over the instances, of the gap of the incumbent solution to the best solution found for that instance in the ten runs, taken as the mean of the ten runs and as the best of them.}
\label{fig6}
\end{figure}

The ILS-BCD matched or improved the GTSP tour on all 235 instances, finding a strictly shorter tour on 234 of them and a tour of the same length on the remaining one. The best of the ten runs is on average $0.441\%$ shorter than $f_{\mathrm{GTSP}}$, with a median of $0.245\%$ and a largest improvement of $5.798\%$. These statistics leave out the six instances whose neighborhoods share a common point, on which the ILS-BCD finds a tour of length zero, yielding a gap of -100\% as defined above. Using the mean across the ten runs instead of the best value preserves virtually the same comparative performance, $0.438\%$ against $0.441\%$ on average, because on 201 of the 235 instances every one of the ten runs reaches the best solution found; of the remaining 34, 18 have more than 50 neighborhoods. This is also why the coefficient of variation in Table~\ref{tab5} rounds to $0.0000\%$ for $m \le 50$. Figure~\ref{fig5} illustrates the solutions obtained by the ILS-BCD for four representative instances of increasing size. Figure~\ref{fig6} shows the progress of the incumbent objective value over time, measured relative to the best solution found across the ten runs for each instance and averaged over all 235 instances. The curves exhibit little improvement after approximately 100 seconds, suggesting that shorter time limits may be sufficient for most instances. Indeed, reducing the time limit to 100 seconds would leave the best-of-ten solution unchanged for 231 instances. For the remaining four instances, all with $m=100$, the best solution would be only $0.061\%$ longer on average.

The 235 instances include the 87 of Section~\ref{sec52}, and for 80 of them the optimal tour length is known, certified by at least one of the two MIQCP formulations, so on this subset the ILS-BCD can be measured against a certified reference rather than against the GTSP tour. On four of the 80 the neighborhoods share a common point and the optimal tour has length zero; the ILS-BCD returns such a tour in all ten runs. On the remaining 76, which have $8 \le m \le 15$, the best of the ten runs differs from the optimal tour length reported in Table~S3 by less than $0.001\%$ on every instance, by $0.000040\%$ on average and by at most $0.000574\%$; the ten runs return the same tour length on 74 of the 76, and differ by less than $4\times10^{-6}\%$ on the other two. On seven instances the tour returned by the ILS-BCD is marginally shorter than the one certified by CPLEX, a difference that lies within CPLEX’s default relative gap tolerance of $0.01\%$, the same effect discussed in Section~\ref{sec53}. Taking that tolerance into account, the lower bounds of Table~S3 certify every tour found by the ILS-BCD to be within $0.011\%$ of optimal, and no tour falls below a lower bound. The computational effort differs by orders of magnitude: every ILS-BCD run stopped by stagnation, taking $6.67$ seconds of CPU time on average and never more than $17.89$, whereas the formulation that proved optimality took $199.23$ seconds on average, and up to $3{,}481.37$.

What the search itself contributes is measured by the improvement of the ILS-BCD over its own initial solution, the multi-start nearest neighbor construction followed by a full BCD of Section~\ref{sec4}, which involves no randomness and is therefore the same in all ten runs. Table~\ref{tab5} reports this improvement beside the gap to the GTSP tour. It is $8.99\%$ on average over the 229 instances whose optimal tour has positive length, and grows monotonically with the number of neighborhoods, from $1.31\%$ for $m \le 10$ to $18.12\%$ for $51 \le m \le 100$: on the smallest instances the construction is already at or near the optimum, and on 29 of them, all with $m \le 25$, the search returns it unchanged. The initial solution is itself $10.361\%$ longer than the GTSP tour on average, so the $0.441\%$ by which the ILS-BCD improves on that tour is the net outcome of a much larger movement.

On average across all runs, the algorithm took 213.8 seconds and reached its best solution after 22.7 seconds. Of the 2,350 runs, 2,132 stopped after 3,000 consecutive iterations without improvement, whereas the remaining 218 reached the time limit—all occurring on instances with $m \geq 49$ and including every run of the thirteen instances with $m=100$. The former completed 3,020 ILS-BCD iterations on average and reached their best solution at iteration 22; thus, most of the execution confirmed a solution found early in the search, in agreement with Figure~\ref{fig6}. The latter reached their best solution at iteration 267 on average, after 230.4 seconds of the 1,200 available.

Table~\ref{tab6} reports how the $286.2$ million neighborhood searches recorded in these runs are distributed among the eight neighborhoods of the VND, with the rows in the order in which the VND explores them. Since the descent returns to the first neighborhood after every improving move, a neighborhood is reached only when all those before it fail on the current solution. Consequently, the number of searches decreases along the VND order, with the first three neighborhoods accounting for $69.1\%$ of all searches. For each neighborhood, the table distinguishes searches that improve the current solution, which we call local improvements, from those that also improve the best solution found so far, and reports the corresponding counts as percentages of the searches performed by that neighborhood and of the respective totals. The two counts differ by nearly three orders of magnitude: $83.1$ million local improvements versus $44{,}258$ best improvements. This large difference is explained by the fact that the best solution is reached at iteration $22$ on average, while each run continues for another $3{,}000$ iterations. Thus, most searches are performed after the best solution found during the run has already been reached. The resulting rates of best improvement are consequently small for all eight neighborhoods, but they still reveal substantial differences in their relative contributions. The 2-Opt stands out in particular: it improves the current solution in $63.5\%$ of its searches and accounts for $59.2\%$ of all best improvements, despite representing only $15.1\%$ of all searches. Insertion follows, improving the current solution in $46.4\%$ of its searches and accounting for $13.7\%$ of all best improvements from only $5.2\%$ of the searches. String relocation shows a similar pattern on a smaller scale, accounting for $6.1\%$ of the best improvements from $2.8\%$ of all searches. Balas--Simonetti and double bridge, the first two neighborhoods in the VND order, account for $54.6\%$ of all local improvements and $17.4\%$ of the best improvements. Thus, they contribute substantially to the descent, even though a relatively small share of the new best solutions comes from these neighborhoods. Swap, 3-Opt, and Gutin have the lowest local-improvement rates, at $6.1\%$, $2.9\%$, and $2.1\%$, respectively, and together account for only $3.7\%$ of the best improvements. Nevertheless, all eight neighborhoods produce best improvements at some point.

\begin{table}[ht!]
\centering
\small
\renewcommand{\arraystretch}{1.15}
\begin{tabular}{lrr@{\hspace{15pt}}rr@{\hspace{15pt}}rrr}
\toprule
\multirow{2}{*}{\raisebox{-0.10ex}{Neighborhood}}
& \multicolumn{2}{c}{Searches}
& \multicolumn{2}{c}{Local improvements}
& \multicolumn{3}{c}{Best improvements} \\
\cmidrule(lr){2-3}
\cmidrule(lr){4-5}
\cmidrule(lr){6-8}
& \raisebox{0.9ex}{Number}
& \shortstack{\% of\\all}
& \raisebox{0.9ex}{Number}
& \shortstack{\% of its\\searches}
& \raisebox{0.9ex}{Number}
& \shortstack{\% of its\\searches}
& \shortstack{\% of all\\improvements} \\
\midrule
Balas--Simonetti
& 89{,}601{,}286 & 31.3
& 25{,}838{,}643 & 28.8
& 4{,}532 & 0.0051 & 10.2 \\

Double bridge
& 63{,}762{,}643 & 22.3
& 19{,}534{,}357 & 30.6
& 3{,}180 & 0.0050 & 7.2 \\

Gutin
& 44{,}228{,}286 & 15.5
& 917{,}795 & 2.1
& 206 & 0.0005 & 0.5 \\

2-Opt
& 43{,}310{,}491 & 15.1
& 27{,}506{,}690 & 63.5
& 26{,}183 & 0.0605 & 59.2 \\

Swap
& 15{,}803{,}801 & 5.5
& 962{,}756 & 6.1
& 1{,}132 & 0.0072 & 2.6 \\

Insertion
& 14{,}841{,}045 & 5.2
& 6{,}888{,}471 & 46.4
& 6{,}062 & 0.0408 & 13.7 \\

String relocation
& 7{,}952{,}574 & 2.8
& 1{,}261{,}349 & 15.9
& 2{,}683 & 0.0337 & 6.1 \\

3-Opt
& 6{,}691{,}225 & 2.3
& 191{,}477 & 2.9
& 280 & 0.0042 & 0.6 \\
\bottomrule
\end{tabular}
\caption{Contribution of each neighborhood over the $2{,}350$ runs, in the order in which the VND explores them. A local improvement is a search that returns a solution better than the current one, whereas a best improvement also improves the best solution found so far. Overall, the $286{,}191{,}351$ searches produced $83{,}101{,}538$ local improvements and $44{,}258$ best improvements. For the 218 runs stopped by the time limit, the counts cover the search up to the last improvement of the best solution found.}
\label{tab6}
\end{table}

\begin{figure}[ht!]
\begin{center}
\begin{tabular}{ccc}
\includegraphics{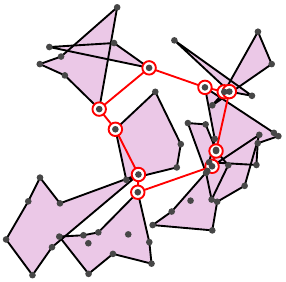} &
\includegraphics{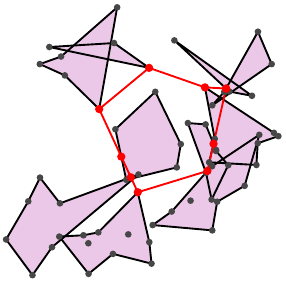} &
\includegraphics{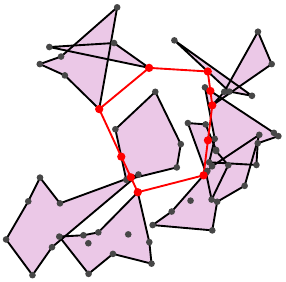} \\ [2mm]
(a) GTSP, $1{,}797.63$ & (b) GTSP-BCD, $1{,}777.16$ & (c) ILS-BCD, $1{,}748.97$
\end{tabular}
\end{center}
\caption{Comparison between the GTSP, the GTSP-BCD and the ILS-BCD solutions on the TSPN instance 10i60-21, which has ten neighborhoods. In (a) the GTSP tour, whose visit points, drawn as white circles, are points of the original GTSP clusters, shown as gray dots. In (b) the same sequence with the visit points recomputed by the BCD, so that they may lie anywhere within the corresponding polygonal neighborhood. In (c) the best solution found by the ILS-BCD, which visits the neighborhoods in a different order.}
\label{fig7}
\end{figure}

The tours found by the ILS-BCD are on average shorter than the GTSP tour, and it is worth asking whether this comes only from moving the visit points inside their neighborhoods, along the sequence of the GTSP tour, or also from visiting the polygonal neighborhoods in a different order. On 122 of the 235 instances the best solution found by the ILS-BCD visits them in the same order as the GTSP tour, and there the gain comes entirely from the optimization of the visit points: the tour is $0.33\%$ shorter than the GTSP tour on average and coincides with the one obtained by keeping the sequence of the GTSP tour and optimizing only its visit points with the BCD. On the remaining 113 instances the ILS-BCD visits the neighborhoods in a different sequence, and the gain is larger: the tour is $0.56\%$ shorter than the GTSP tour on average, and it is also $0.10\%$ shorter than the one that the BCD obtains from the sequence of the GTSP tour, which it improves on 83 of them and matches on the remaining 30. On these the sequences differ but the two tours describe the same closed curve, with some of the visit points falling at different places along it. Figure~\ref{fig7} shows one instance where the solution tour given by the ILS-BCD visits the neighborhoods in a different sequence, and the tour is $2.71\%$ shorter than the GTSP tour and $1.59\%$ shorter than the one that the BCD obtains from the sequence of the GTSP tour.

\begin{figure}[ht!]
\begin{center}
\begin{tabular}{ccc}
\includegraphics{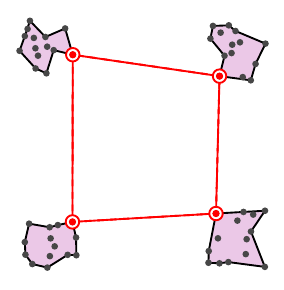} &
\includegraphics{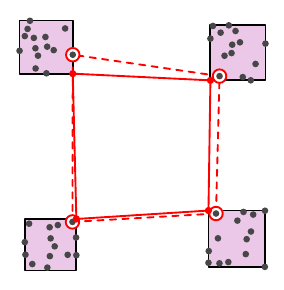} &
\includegraphics{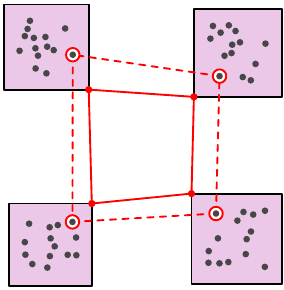} \\ [2mm]
(a) gap $0.0\%$ & (b) gap $-8.54\%$ & (c) gap $-30.24\%$ \\ [4mm]
\includegraphics{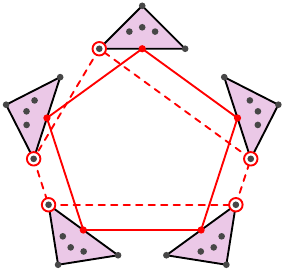} &
\includegraphics{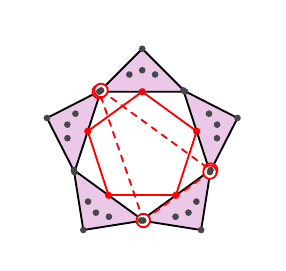} &
\includegraphics{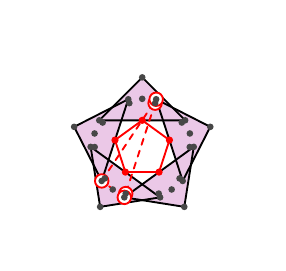} \\ [2mm]
(d) gap $-1.73\%$ & (e) gap $-6.07\%$ & (f) gap $-26.59\%$
\end{tabular}
\end{center}
\caption{The gain over the GTSP tour increases as the neighborhoods become larger relative to the distances between them. In the first row, the neighborhoods are enlarged while the clusters remain fixed: (a) shows the neighborhoods obtained from the conversion described in Section~\ref{sec51}, for which the two tours coincide; (b) replaces each neighborhood with its smallest enclosing square; and (c) enlarges these squares by a factor of $1.6$. Since the cluster points remain fixed, the optimal GTSP tour is the same in all three panels. In the second row, the neighborhoods remain fixed while the five clusters are moved towards the center, with their original distances to the center scaled by factors of $1$, $0.7$, and $0.5$ in (d), (e), and (f), respectively. The dashed lines represent the optimal GTSP tour, and the solid lines represent the best solutions found by ILS-BCD. The corresponding gaps relative to the GTSP tour are reported below each panel.}
\label{fig8}
\end{figure}

The gains reported above are bounded by the way the TSPN instances were derived in Section~\ref{sec51}. Each neighborhood is the convex hull of a GTSP cluster with at most three of the discarded points reinserted as reflex vertices, so it barely extends beyond the points it was built from. Neighborhoods generated with a larger extent would therefore have produced larger gains. The effect is easy to isolate, because enlarging a neighborhood that already contains its cluster keeps every solution of the GTSP instance feasible and leaves $f_{\mathrm{GTSP}}$ unchanged. Figure~\ref{fig8} shows it on an instance with four clusters built by the procedure of Section~\ref{sec51}: with the converted neighborhoods the ILS-BCD returns the GTSP tour itself, a gap of $0.0\%$ (Figure~\ref{fig8}(a)); replacing each neighborhood by the smallest square that contains it takes the gap to $-8.54\%$ (Figure~\ref{fig8}(b)), and enlarging that square further takes it to $-30.24\%$ (Figure~\ref{fig8}(c)). The same effect can be observed by keeping the neighborhoods fixed while moving the clusters closer to one another, thereby also bringing the polygons derived from them closer together. This is because the gain depends on the size of the neighborhoods relative to the distances between them. The second row of Figure~\ref{fig8} illustrates this effect for an instance with five clusters. Specifically, scaling the clusters' distances from the center to $70\%$ and then to $50\%$ of their original values improves the gain from $-1.73\%$ in Figure~\ref{fig8}(d) to $-6.07\%$ in Figure~\ref{fig8}(e) and $-26.59\%$ in Figure~\ref{fig8}(f).

\begin{figure}[ht!]
\begin{center}
\begin{tabular}{cc}
\includegraphics[width=0.46\textwidth]{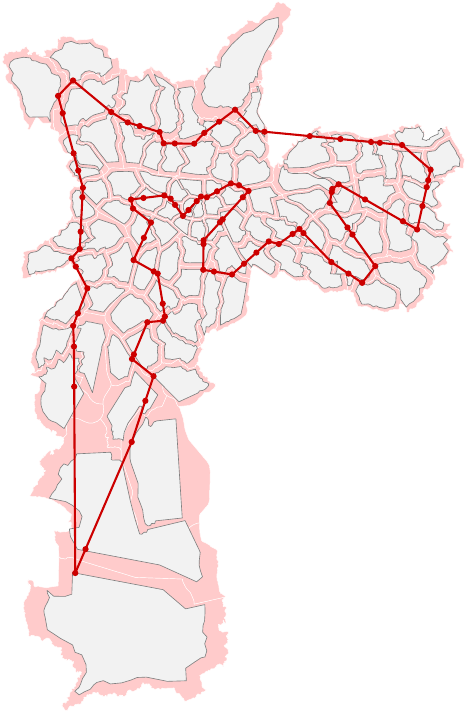} &
\includegraphics[width=0.46\textwidth]{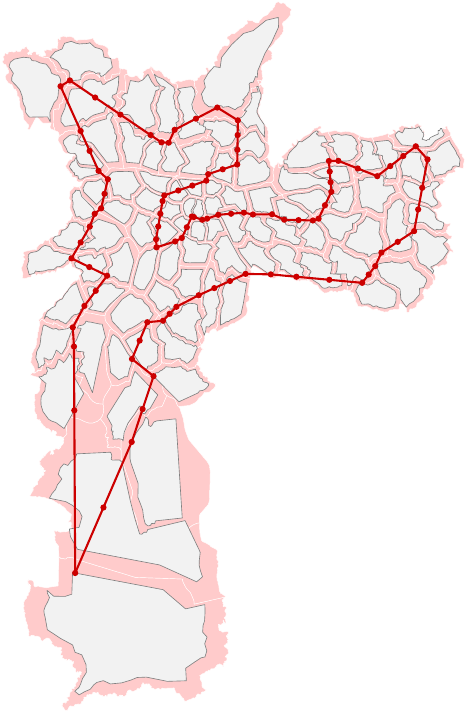} \\ [2mm]
(a) Solution of~\cite{Birgin2022}, $212{,}292.01$ & (b) ILS-BCD, $199{,}362.73$
\end{tabular}
\end{center}
\caption{The instance of the city of São Paulo introduced in~\cite{Birgin2022}. In (a), the solution reported there, of length $212{,}292.01$; in (b), the best solution found by the ILS-BCD from the same initial solution, of length $199{,}362.73$. The 96 polygonal neighborhoods to be visited are shown in gray and the polygons of the original map are shown in coral in the background.}
\label{fig9}
\end{figure}

We close this section with an instance of a different nature. In~\cite{Birgin2022}, the authors introduce a TSPN instance built from the map of the city of São Paulo, Brazil. This instance consists of 96 polygonal neighborhoods with between 51 and 57 vertices each, 4{,}966 in total. They are far less regular than the polygonal neighborhoods converted in Section~\ref{sec51}, with 16 to 29 reflex vertices each instead of at most three, and all but two pairs of them are disjoint. The tour reported there has length $212{,}292.01$ and was obtained with $3{,}589.31$ seconds of CPU time. Running the code of~\cite{Birgin2022} on the machine used here reproduces that tour in about a third of that time, so the CPU time limit of $1{,}200$ seconds of our protocol is comparable to the time reported there. We ran the ILS-BCD from the same initial solution, of length $229{,}139.65$, given by the constructive heuristic of~\cite{Birgin2022}. Eight of the ten runs returned a tour between $199{,}362.73$ and $199{,}362.78$, $6.090\%$ shorter than the length reported there, and the other two returned $199{,}390.36$ and $199{,}637.10$, $6.077\%$ and $5.961\%$ shorter. Figure~\ref{fig9} shows the solution of~\cite{Birgin2022} and the best solution found by the ILS-BCD.

\subsection{Larger instances and the limits of the ILS-BCD} \label{sec56}

\begin{figure}[ht!]
\centering
\includegraphics{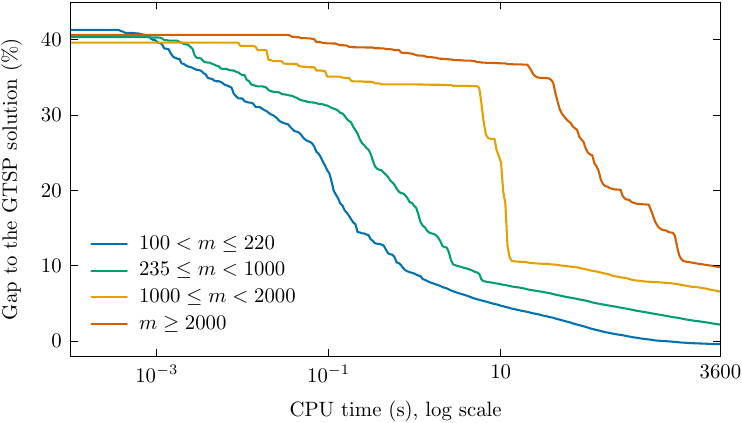}
\caption{Mean gap to the GTSP reference tour over CPU time, averaged over the ten runs of each instance, for the four size groups of Table~\ref{tab7}.} 
\label{fig10}
\end{figure}

The instances of Section~\ref{sec51} were restricted to at most 100 neighborhoods. The conversion procedure applies unchanged to the larger instances of the two libraries. Thus, we converted the remaining 60 instances with more than 100 clusters of the GTSP libraries considered in Section~\ref{sec51} for which a reference tour is available, 45 from the GTSP Lib and 15 from the MOM Lib, leaving aside one with $m = 14{,}202$. On these instances the number of polygons varies between 115 and 6{,}325 and the total number of vertices between 457 and 30{,}939. Each instance was solved in ten independent runs with the configuration of Section~\ref{sec54}, the same used throughout, and a CPU time limit of one hour. Table~S9 in the supplementary material reports the results for individual instances, Table~\ref{tab7} summarizes them by instance size, and Figure~\ref{fig10} shows the evolution of the incumbent solution over time.

\begin{table}[ht!]
\centering
\small
\setlength{\tabcolsep}{4pt}
\renewcommand{\arraystretch}{1.15}
\resizebox{\textwidth}{!}{
\begin{tabular}{lrrrrrrr}
\toprule
\multirow{2}{*}{Instance size ($m$)}
& \multirow{2}{*}{\raisebox{-0.10ex}{\# Inst.}}
& \multicolumn{2}{c}{Instances with shorter tour than}
& \multicolumn{3}{c}{Gap to the GTSP tour (\%)}
& \multirow{2}{*}{\shortstack[r]{Mean \# ILS\\iterations}} \\
\cmidrule(lr){3-4} \cmidrule(lr){5-7}
& & \multicolumn{1}{c}{GTSP}
& \multicolumn{1}{c}{GTSP-BCD}
& \multicolumn{1}{r}{Best}
& \multicolumn{1}{r}{Mean}
& \multicolumn{1}{r}{Worst} & \\
\midrule
$100 < m \le 220$   & 26 & \multicolumn{1}{c}{25} & \multicolumn{1}{c}{24} & $-0.575$ & $-0.404$ & $-0.098$ & 1{,}093 \\
$235 \le m < 1000$  & 22 &  \multicolumn{1}{c}{3} &  \multicolumn{1}{c}{2} & $+1.430$ & $+2.195$ & $+2.989$ &    382 \\
$1000 \le m < 2000$ &  2 &  \multicolumn{1}{c}{0} &  \multicolumn{1}{c}{0} & $+5.625$ & $+6.538$ & $+7.380$ &     55 \\
$m \ge 2000$        & 10 &  \multicolumn{1}{c}{0} &  \multicolumn{1}{c}{0} & $+9.634$ & $+9.805$ & $+9.970$ &      0 \\
\bottomrule
\end{tabular}}
\caption{ILS-BCD results on the 60 instances with more than 100 neighborhoods, grouped by instance size. Each run has a CPU time limit of one hour. The two counts report the number of instances for which the best of the ten runs yields a tour strictly shorter than the GTSP tour and the GTSP-BCD tour, respectively. The gap columns report the mean, over the instances in each group, of the gaps corresponding to the best, mean and worst tour lengths across the ten runs. The last column reports the mean number of ILS iterations completed per run up to the last improvement of the best solution found.} 
\label{tab7}
\end{table}

Up to about 220 neighborhoods the method behaves as it does in Section~\ref{sec55}. It returns a shorter tour than the GTSP reference on 25 of the 26 instances, and on 24 of them it also improves on the GTSP-BCD tour. This is not an artifact of the best of the ten runs: averaged over the group, even the worst of them is shorter than the reference. With the time limit of 1{,}200 seconds used in Section~\ref{sec55} the same group still improves on the GTSP tour on 24 of the 26 instances, but the worst of the ten runs is then $0.368\%$ longer than the reference rather than $0.098\%$ shorter, so the longer limit improves the consistency of the ten runs rather than the quality of the best among them. Beyond that size the reference is shorter, by a margin that grows with $m$. On the 22 instances with $235 \le m < 1000$ the ILS-BCD is shorter than the GTSP tour on only three, and by at most $0.175\%$. On the 12 instances with $m \ge 1000$ it is never shorter. The last column of Table~\ref{tab7} identifies what happens. The number of ILS-BCD iterations completed before the last improvement falls from more than a thousand to a few dozen and then to none: on the ten instances with $m \ge 2000$, all ten runs of all ten instances record a single call to the VND and no ILS-BCD iteration at all. The hour is spent on the construction of the initial solution and on one descent of the local search that does not finish. The construction alone, which is a nearest neighbor tour started from each of the $m$ neighborhoods, each followed by a complete BCD, takes 32 seconds at $m = 2{,}000$ and 990 seconds at $m = 6{,}325$, the latter 27\% of the budget. What these instances measure is therefore the construction followed by one incomplete local search, not the metaheuristic.

The comparison with the GTSP literature locates the difficulty. The reference tours used here were produced by a specialized GTSP solver on instances with up to 6{,}325 clusters~\cite{Helsgaun2015}, and instances of this size are standard benchmarks in that literature~\cite{Schmidt2022}: when each neighborhood is a finite set of points, thousands of them are within reach. What our method does not absorb at this scale is the replacement of that finite set by a non-convex polygon. Every candidate sequence then requires the visit points to be recomputed, and although the BCD performs that computation in microseconds when $m$ is a hundred, the cost of a single pass of the VND over its eight neighborhoods grows until it exceeds the budget. The challenge these instances pose can therefore be stated precisely: to devise a method that, at this scale, does better than treating the TSPN as a GTSP and applying the BCD once to the resulting sequence---which is exactly the reference that the ILS-BCD does not beat here.

One direction seems worth exploring, and it attacks the size rather than the constant. The heuristics of~\cite{Mennell2011} and~\cite{Wang2019} for the close-enough TSP are built on Steiner zones, the non-empty intersections of several disks: a point in a zone of degree $k$ visits $k$ cities at once, so identifying the zones and choosing one point in each turns the problem into a Traveling Salesman Problem on the chosen points. Transposed to polygonal neighborhoods, the construction would reduce the number of entities to be sequenced wherever neighborhoods overlap, and it is that number, not the shape of each neighborhood, that exhausts the budget in Table~\ref{tab7}. Two properties of our setting make the transposition natural. The neighborhoods of Section~\ref{sec2} are unions of convex polygons, and the intersection of two such sets is the union of the intersections of their convex pieces, so the class is closed under intersection and a zone is again an object that the formulation and the BCD already handle. And the zones reduce the overlap that invalidates the second hypothesis of Section~\ref{sec3}, since neighborhoods that meet can be merged into a single zone. Whether the number of zones can be kept manageable, how to choose among zones of different degrees, and whether the sequence they induce is good enough to be refined rather than rebuilt are open questions, and the answer to none of them is obvious. What makes the direction attractive is that the expensive half of the zone-based methods---recomputing the visit points by second-order cone programming---is the half that the BCD of Section~\ref{sec3} performs several orders of magnitude faster.

\section{Conclusions and future work} \label{sec6}

In this work, we addressed the Traveling Salesman Problem with Neighborhoods (TSPN) in the Euclidean plane, specifically focusing on general non-convex polygonal neighborhoods. To tackle the continuous and combinatorial challenges inherent to this setting, we presented a Mixed-Integer Quadratically Constrained Programming (MIQCP) formulation along with an efficient hybrid metaheuristic framework.

The proposed MIQCP model handles complex polygonal regions by decomposing them into unions of convex subpolygons. By incorporating linear bounding inequalities derived from vector norm relations, we reduced exact solver CPU resolution times by approximately 22\% on instances solved to optimality. To address the continuous subproblem of determining optimal visit points along a fixed tour permutation, we introduced a Block Coordinate Descent (BCD) procedure that resolves fast one-dimensional geometric subproblems. Empirical evaluations demonstrated that BCD achieves visit points of similar quality to commercial solvers while running more than six orders of magnitude faster. This gain in computational efficiency proved essential for enabling the hybrid Iterated Local Search (ILS) metaheuristic—accelerated via candidate lists and truncated BCD evaluations within a Variable Neighborhood Descent (VND) scheme to thoroughly explore the combinatorial search space at a low computational cost. Tested on 235 benchmark instances adapted from the GTSP literature, the ILS-BCD algorithm matched or improved upon best-known GTSP tour lengths in all cases, obtaining strictly shorter tours in 234 instances. On the 80 instances whose optimal tour length is known from the MIQCP formulations, it attained that length in every one, within a relative difference below $0.001\%$, and did so in $6.67$ seconds of CPU time on average against $199.23$ for the formulation that proved optimality. On the larger instances of the same libraries, with up to 6{,}325 neighborhoods, the method remains competitive up to about 220 neighborhoods and is outperformed by the GTSP reference beyond that size, where a single pass of the local search no longer fits within an hour.

Future research offers several promising directions. Extending the BCD scheme to multi-point blocks comprising two or more consecutive visit points could yield higher solution quality while balancing the additional subproblem complexity. Additionally, adapting the ILS-BCD framework to three-dimensional TSPN variants or dynamic settings involving moving polygonal neighborhoods represents an intriguing avenue. A further direction, suggested by Section~\ref{sec56}, is to adapt to polygonal neighborhoods the Steiner-zone constructions used for the close-enough TSP, which reduce the number of entities to be sequenced by merging neighborhoods that overlap. Finally, integrating BCD-generated upper bounds and valid inequalities into exact branch-and-bound or branch-and-cut algorithms could further advance exact solution methods for the TSPN.\\

\noindent
\textbf{Acknowledgment:} We sincerely thank Marcelo Marazzi for his careful editing and helpful suggestions that enhanced the quality of the written text in English.

\bibliographystyle{plain}
\bibliography{ber2025}

\end{document}